\documentstyle[leqno]{article}
\newcommand{\se}[1]{{\section{#1}} {\setcounter{equation}{0}}}
\newtheorem{th}{Theorem}[section]
\newtheorem{lm}{Lemma}[section]
\newtheorem{prop}{Proposition}[section]

\def\k{{K\"{a}hler }}

\def\cy{{Calabi-Yau }}
\def\l{{Lagrangian }}
\def\sl{{special Lagrangian }}
\input epsf
\begin{document}
\hbadness=10000
\title{{\bf Lagrangian torus fibration of quintic \cy hypersurfaces I:}\\
{\Large {\bf Fermat quintic case}}}
\author{Wei-Dong Ruan\\
Department of mathematics\\
Columbia University\\
New York, NY 10027\\
ruan@math.columbia.edu}
\date{Revised December 1999}
\footnotetext{Partially supported by NSF Grant DMS-9703870.}
\maketitle
\begin{abstract}
In this paper we give a construction of Lagrangian torus fibration for Fermat type quintic \cy hypersurfaces via the method of gradient flow. We also compute the monodromy of the expected \sl torus fibration and discuss structures of singular fibers.
\end{abstract}
\tableofcontents
\se{Introduction and background}
In this paper we give a construction of Lagrangian torus fibration for Fermat type quintic \cy manifolds via the method of gradient flow. This method will produce Lagrangian torus fibration for \cy hypersurfaces in general toric varieties. The results for general quintic hypersurfaces appeared in \cite{lag3}; for \cy hypersurfaces in general toric varieties are written in \cite{tor}.\\\\
The motivation of our work comes from the study of Mirror Symmetry. Mirror Symmetry conjecture originated from physicists' work in conformal field theory and string theory. It proposes that for a \cy 3-fold $M$ there exists a \cy 3-fold $W$ as its mirror manifold. The quantum geometry of $M$ and $W$ are closely related. In particular one can compute the number of rational curves in $M$ by solving the Picard-Fuchs equation coming from variation of Hodge structure of $W$. \\\\
The first example of mirror symmetry was the Fermat type quintic given in the calculations by Candelas et al (\cite{Can}). Numerous examples of such calculations were worked out after that. A more general construction via toric varieties was given by Batyrev (\cite{Bat}). For a more complete history and reference the readers can consult \cite{mirrorbook}.\\\\
In 1996 Strominger, Yau and Zaslow (\cite{SYZ}) proposed a geometric construction of mirror manifold via special Lagrangian torus fibration. According to their program (we will call it SYZ construction), a \cy 3-fold should admit a special Lagrangian torus fibration. The mirror manifold can be obtained by dualizing the fibers. Or equivalently, the mirror manifold of $M$ is the moduli space of special Lagrangian 3-torus in $M$ with a flat $U(1)$ connection. This conjectural construction was the first to give the mirror manifold directly from a \cy manifold itself. \\\\
The notion of special Lagrangian submanifolds was first given by Harvey and Lawson in their celebrated paper \cite{HL}. For a \cy manifold $(X, \Omega, \omega_g)$, where $\Omega$ is a holomorphic $(n,0)$ form and $\omega_g$ is the \k form of the \cy metric, a \l submanifold $L$ in $X$ is called \sl if $L$ is area minimizing, or equivalently $\Omega$ restrict to $L$ is a constant multiple of the volume form of $L$.\\\\
According to the SYZ construction, \sl submanifold and \sl fibration for \cy manifolds seem to play very important roles in understanding mirror symmetry. However, despite its great potential in solving the mirror symmetry conjecture, our understanding in special Lagrangian submanifolds is very limited. The known examples are mostly explicit local examples or examples coming from $n=2$. There are very few examples of \sl submanifold or \sl fibration for dimension higher than two. M. Gross, P.M.H. Wilson, N. Hitchin, P. Lu and R. Bryant(\cite{Gross1}\cite{Gross2}\cite{GW}\cite{H}) did some work in this area in recent years. On the other extreme, in \cite{Z}, Zharkov constructed some non-\l torus fibration of \cy hypersurfaces in toric variety.\\\\
When dimension $n=2$, \cy manifold is hyper\k. Therefore for any \cy metric, there is an $S^2$ family of complex structures that are compatible with the given \cy metric. Special Lagrangian submanifold for one compatible complex structure exactly corresponds to complex curves for another compatible complex structure. Therefore \sl theory is reduced to the theory of complex curves in $X$, which is fairly well understood. For $n\geq 3$, there are no nice interpretations like this. \\\\
Given our lack of knowledge for \sl, one may consider relaxing the requirement to consider \l fibration. Special Lagrangians are very rigid and hard to find. On the other hand, \l submanifolds are more flexible and can be modified locally by Hamiltonian deformation. This is a reasonable first step to take. For many applications to mirror symmetry, especially those concerning (symplectic) topological structure of fibration, \l fibration will provide quite sufficient information. In this paper we mainly consider Lagrangian torus fibration of \cy hypersurfaces in toric varieties. \\\\
Our idea is a very natural one. We try to use gradient flow to get Lagrangian torus fibration from a known Lagrangian torus fibration at the "Large Complex Limit". It will in principle be able to produce \l fibration in general \cy hypersurfaces. For simplicity, it is helpful to explore the case of Fermat type quintic \cy threefold in ${\bf CP^4}$ in detail first, which is the case studied in \cite{Can}. Most of the essential features for \cy hypersurfaces in general toric varieties already show up here.\\\\
The paper is organized as follows: In Section 2 we will first describe a \l torus fibration for  Fermat type quintic \cy familly $\{X_\psi\}$ in ${\bf P^4}$\\
\[
p_{\psi}=\sum_1^5 z_k^5 - 5\psi \prod_{k=1}^5 z_k=0
\]\\
at the large complex limit $X_{\infty}$\\
\[
p_{\infty}=\prod_{k=1}^5 z_k=0.
\]\\
In Section 3 we construct an appropriate vector field whose gradient flow will give \l fibration to nearby \cy hypersurfaces. In Section 4, we will discuss the expected fibration structure by explicitly computing the monodromy transformations of the expected fibration if the fibration is actually \sl fibration. In Section 5 we will  describe the expected structures of singular fibers implied by monodromy information. In Section 6 we will compare our \l fibration with the expected \sl fibration---their topological structure are not exactly the same; we will discuss the cause for their differences. Finally in Section 7, we will discuss the relevance to mirror construction for toric \cy manifolds through the dual polyhedra construction.\\\\
This paper is the first part of several papers in this subject. In \cite{lag3} we discuss general construction of  Lagrangian torus fibration for general quintic (non-Fermat) \cy hypersurfaces.  In \cite{lag2} we address the technical problem of modifying our construction into a \l torus fibration with the topological type of the expected one. \\\\

\se{\l torus fibration for large complex limit}
Consider the well studied  Fermat type quintic \cy familly $\{X_\psi\}$ in ${\bf P^4}$ defined by\\
\[
p_{\psi}=\sum_1^5 z_k^5 - 5\psi \prod_{k=1}^5 z_k=0.
\]\\
When $\psi$ approaches $\infty$, The familly approach its ``large complex limit'' $X_{\infty}$\\
\[
p_{\infty}=\prod_{k=1}^5 z_k=0.
\]\\
$X_{\infty}$ is a union of five ${\bf P^3}$. There is a natural degenerate  $T^3$ fibration structure for $X_{\infty}$. Let $\{P_i\}_{i=1}^5$ be five points in ${\bf R^4}$ that are in general position. Consider the natural map $F: \ P^4\longrightarrow R^4$.\\
\[
F([z])=\sum_{k=1}^5 \frac{|z_k|^2}{\sum_{i=1}^5 |z_i|^2}P_k.
\]\\
Then $\Delta={\rm Image}(F)$ is a 4-simplex. $X_{\infty}$ is naturally fibered over $\partial \Delta$ via this map $F$ with general fiber being $T^3$. This is precisely the $T^3$ special Lagrangian fibration for $X_{\infty}$ as indicated by SYZ construction.\\\\
To see this fact, we need to construct a suitable \cy metric on $X_{\infty}$ such that $F$ is a \sl torus fibration with respect to the \k form of the \cy metric. Since $X_{\infty}$ is just a union of several ${\bf CP^3}$'s, we will concentrate on one of these ${\bf CP^3}$. We will actually carry out the discussion for general ${\bf CP^n}$.\\\\
In general, on $ {\bf C^{n+1}} $ there is a natural flat \k metric with \k form\\
\[
\omega_0 = \sum_{i=1}^{n+1} \frac{dz_i\wedge d\bar{z}_i}{|z_i|^2}
\]\\
which is singular along the coordinate hyperplanes. Consider the natural projection $\pi: \ {\bf C^{n+1}}-\{0\}\rightarrow {\bf CP^{n}} $. The restriction of $\omega_0$ to the hypersurface $\prod_{i=1}^{n+1}z_i =1$ in $\bf C^{n+1}$ naturally push forward via $\pi$ to a flat \k metric $\omega$ on ${\bf CP^{n}}$ that is singular along the union of projective coordinate hyperplanes, which is exactly our large complex limit $Y_{\infty}$. (Here we use $Y_{\infty}$ to denote dimention $n-1$ version of the large complex limit in ${\bf CP^{n}}$ to distinguish from $X_\infty$ that corresponding to $n=4$ case.) Take local coordinate $x_i=\frac{z_i}{z_{n+1}}$ for $i=1,\cdots, n$. Then\\
\[
\omega = \sum_{i=1}^{n} \frac{dx_i\wedge d\bar{x}_i}{|x_i|^2} - \frac{1}{n+1} \left(\sum_{i=1}^{n} \frac{dx_i}{x_i}\right)\left(\sum_{i=1}^{n} \frac{d\bar{x}_i}{\bar{x}_i}\right)
\]\\
It is nice to compute their \k potentials:\\
\[
\omega_0 = \partial \bar{\partial} \left( \sum_{i=1}^{n+1} (\log|z_i|^2)^2 \right)
\]\\
\[
\omega = \partial \bar{\partial} \left( \sum_{i=1}^{n} (\log|x_i|^2)^2 - \frac{1}{n+1} (\sum_{i=1}^{n}\log|x_i|^2)^2\right)
\]\\
Clearly\\
\[
\omega^n = (-1)^{\frac{n(n-1)}{2}} \frac{n!}{n+1}\Omega \wedge \bar{\Omega}
\]\\
where\\
\[
\Omega = \bigwedge_{i=1}^n \frac{dx_i}{x_i}
\]\\
is the holomorphic n-form on ${\bf CP^{n}}\backslash Y_{\infty}$. Notice that $\omega$ corresponds to a complete \cy metric on ${\bf CP^{n}}\backslash Y_{\infty}$. Consider the n-dimension real torus $T^n$ represented by\\
\[
e^{i\theta} = (e^{i\theta_1}, e^{i\theta_2}, \cdots, e^{i\theta_n})
\]\\
as a real abelian Lie group. Then $T^n$ act naturally on ${\bf CP^{n}}$\\
\[
e^{i\theta}(x) = (e^{i\theta_1}x_1, e^{i\theta_2}x_2, \cdots, e^{i\theta_n}x_n)
\]\\
as symplectomorphisms. The moment map is exactly $F: {\bf CP^{n}}\backslash Y_{\infty} \rightarrow {\bf R^n}$\\
\[
F(x)= ( \log |x_1|^2, \log |x_2|^2, \cdots, \log |x_n|^2)
\]\\
It is easy to see that when $\Omega$ is restricted to the fibre, we have\\
\[
\Omega|_{T^n(x)} = i^n\bigwedge_{i=1}^n d\theta_i
\]\\
Therefore $F$ naturally gives us the \sl fibration of ${\bf CP^{n}}\backslash Y_{\infty}$. Specialize to $n=3$ case, this construction gives us the \sl torus fibration with respect to a complete (flat) \cy metric on $X_\infty\backslash {\rm Sing}(X_\infty)$.\\\\
When $\psi$ is large, according to SYZ conjecture, we {\it expect} that $X_{\psi}$ will also possess a special largrangian $T^3$ fibration with base identified with $\partial \Delta$ which is topologically an $S^3$. There are still serious analysis and geometric works to be done to totally justify the special Lagrangian fibration. We will instead give a natural construction of Lagrangian torus fibration.\\\\
For our purpose, we will consider the Fubini-Study metric\\
\[
\omega_{FS} = \partial \bar{\partial} \log(1+|x|^2).
\]\\
The nice thing about the Fubini-Study metric is that the $T^n$ action is also a symplectomorphism with respect to $\omega_{FS}$. The corresponding moment map is\\
\[
F(x)= \left( \frac{|x_1|^2}{1+|x|^2}, \frac{|x_2|^2}{1+|x|^2}, \cdots, \frac{|x_n|^2}{1+|x|^2}\right)
\]\\
which is easy to see if we write $\omega_{FS}$ in polar coordinates.\\
\[
\omega_{FS} = \partial \bar{\partial} \log(1+|x|^2)=i\sum_{j=1}^n d\theta_i \wedge d\left(\frac{|x_i|^2}{1+|x|^2}\right).
\]\\
The most symmetric expression of the \k potential for the Fubini-Study metric is\\
\[
h= \frac{1 + |x|^2}{\left(\prod_{i=1}^n|x_i|^2\right)^{\frac{1}{n+1}}} = \frac{|z|^2}{\left( \prod_{i=1}^{n+1}|z_i|^2 \right)^{\frac{1}{n+1}}}
\]\\
This expression can also be derived by restricting $|z|^2$ to the hypersurface $\prod_{i=1}^{n+1}z_i=1$ and then push down by $\pi$.\\\\
In general, we have the following lemma:\\
\begin{lm}
For a \k metric\\
\[
\omega = \partial \bar{\partial} h(x)
\]\\
$T^n$ acts as symplectomorphism if and only if $h$ can be chosen to depend only on $|x_i|^2$.
\end{lm}
{\bf Proof:} Assume that this is the case, then\\
\[
\omega = \partial \bar{\partial} h = i\sum_{j=1}^n d\theta_i \wedge d h_i.
\]\\
where\\
\[
h_i=|x_i|^2\frac{\partial h}{\partial |x_i|^2}.
\]\\
From above expression of $\omega$, clearly $T^n$ acts as symplectomorphism and\\
\[
F_h(x) = (h_1, h_2, \cdots, h_n)
\]\\
is the moment map.
\begin{flushright} $\Box$ \end{flushright}
For the purpose of this paper, any $T^n$-invariant \k metric as in Lemma 2.1 is as good. We will mainly use the Fubini-Study metric.\\\\

\se{Lagrangian torus fibration via a gradient flow}
\subsection{The gradient vector field}
Consider the meromorphic function\\
\[
s= \frac{\prod_{k=1}^5 z_k}{\sum_1^5 z_k^5}
\]\\
defined on ${\bf P^4}$. Let $\omega$ denote the \k form of a \k metric $g$ on ${\bf P^4}$, and $\nabla f$ denote the gradient vector field of real function $f=Re(s)$ with respect to the \k metric $g$. To describe the construction, we need the following facts.\\
\begin{lm}
The gradient flow of $f$ leaves the set $\{Im(s)=0\}$ invariant and deforms \l submanifolds in $X_{\infty}$ to \l submanifolds in $X_{\psi}$.\\
\end{lm}
{\bf Proof:} Clearly $\nabla f$ is always perpendicular to level sets of $f$. Therefore,\\
\[
J\nabla f (f)=0.
\]\\
Let $s = f+ ih$. Since $s$ is holomorphic, we have\\
\[
\nabla f (h) = - J\nabla f (f)=0.
\]\\
Therefore $h$ is constant along the gradient flow, or in another word, $\{Im(s)={\rm constant}\}$ is invariant under the gradient flow of $f$. In particular, $\{Im(s)=0\}$ is invariant under the gradient flow of $f$.\\\\
Notice the fact that\\
\[
\{Im(s)=0\} = \bigcup_{\psi ^{-1}\in {\bf R}} X_{\psi}
\]\\
Let $L$ be a \l submanifold of $X_{\infty}$. To prove the lemma, it is sufficient to show that the 4-dimension submanifold $S$ swept out by $L$ under the gradient flow of $f$ is \l submanifold of ${\bf P^4}$. \\\\
First, for a vector field $v$ on $L$, since $L\subset X_{\infty}$ and $\nabla f$ is perpendicular to level sets of $f$ (in particular, perpendicular to $X_{\infty}$), we see that $Jv$ is along $X_{\infty}$ and $\nabla f$ is perpendicular to $Jv$ along $L$. Therefore, $\omega(v,\nabla f)=0$. \\\\
Let $V$, $W$ denote vector fields on $S$ that are invariant under the gradient flow $\phi_t$. $\nabla f$ is a vector field of this type. Using the fact that $f$ is pluri-harmonic and $\omega_g$ is \k form, it is easy to derive that\\
\begin{eqnarray*}
\frac{d}{dt}\omega_g(V,W) &=& -{\cal L}_{\nabla f} \omega_g(V,W) = -(d i(\nabla f)\omega_g)(V,W)\\
&=& (dJdf)(V,W) = -2i\partial \bar{\partial} f(V,W)=0.
\end{eqnarray*}
Therefore $\omega_g(V,W)$ is constant along the flow. Since initial value of invariant vector fields are spanned by vector fields along $L$ and $\nabla f$, by the \l property of $L$ and the fact that $\omega(v,\nabla f)=0$ for vector field $v$ on $L$, we have $\omega_g(V,W)=0$, therefore $S$ is \l.
\begin{flushright} $\Box$ \end{flushright}
{\bf Remark:}\\
(i) The lemma can also be understood roughly by the fact: $\nabla f = H_h$ (the hamiltonian vector field generated by $h=Im(s)$). At the smooth part of the vector field, this fact essentially implies the lemma, although at singular part of the vector field, additional argument as in the proof is needed.\\\\
(ii) The proof of the lemma actually implies something more. Lower dimension torus in $X_{\infty}$ form critical set of $f$. The proof of the lemma also implies that the stable manifold of a lower dimensional torus with respect to flow of $\nabla f$ is \l in ${\bf P^4}$, therefore intersect with $X_{\psi}$ at a \l submanifold.
\begin{flushright} $\Box$ \end{flushright}
With this lemma in mind, the construction of \l torus fibration of $X_{\psi}$ for $\psi$ large is immediate. Recall that in last section we had a canonical Lagrangian torus fibration of $X_{\infty}$ over $\partial \Delta$. Deform along gradient flow of $f$ will naturally induce a \l torus fibration of $X_{\psi}$ over $\partial \Delta$ for $\psi$ large and real.\\\\
However $X_{\infty}$ is singular and $\nabla f$ is also singular where $s$ is singular. To get a really honest \l torus fibration, we need to discuss how to deal with these singularities.\\\\
Along the gradient flow of $f$ we have\\
\[
\frac{df}{dt} = \nabla f \cdot \frac{dx}{dt} = \nabla f \cdot \nabla f = |\nabla f|^2
\]\\
We see the gradient flow of $f$ does not exactly move $X_{\infty}$ to $X_{\psi}$. To ensure this property, the flow has to satisfy $\frac{df}{dt} = {\rm constant}$. This will be true if we scale the vector field $\nabla f$ to\\
\[
V=\frac{\nabla f}{|\nabla f|^2}.
\]\\
\subsection{Flow from the smooth part of $X_{\infty}$}
To understand the flow of $\nabla f$ and $V$, it is helpful to express them in local coordinate. Recall that $X_{\infty}=\cup_{k=1}^5 D_k$, where $D_k=\{z_k=0\}$. Let's choose coordinate $x_i=z_i/z_5$, for $i=1,2,3,4$ and consider $D_4$, where $x_4=0$. Under this coordinate, we have\\
\[
s=\frac{\prod_{i=1}^4x_i}{\sum_{i=1}^4x_i^5 +1}
\]\\
\[
ds =\frac{\prod_{i=1}^4x_i}{\left(\sum_{i=1}^4x_i^5 +1\right)^2}\sum_{i=1}^4 \left(\sum_{j=1}^4x_j^5 +1-5x^5_i\right)\frac{dx_i}{x_i}
\]\\
\[
ds|_{D_4} =\frac{\prod_{i=1}^3x_i}{\sum_{i=1}^3x_i^5 +1}dx_4
\]\\
For simplicity, we choose \k metric to be\\
\[
g=\sum_{i=1}^4 dx_i\wedge d\bar{x}_i
\]\\
Then when restricted to $D_4$, we have\\
\[
|ds|^2 = \left|\frac{\prod_{i=1}^3x_i}{\sum_{i=1}^3x_i^5 +1}\right|^2
\]\\
\[
\nabla s = \frac{\prod_{i=1}^3x_i}{\sum_{i=1}^3x_i^5 +1} \frac{\partial}{\partial \bar{x}_4}
\]\\
These give us\\
\[
V=\frac{\nabla f}{|\nabla f|^2} = Re\left(\frac{\sum_{i=1}^3x_i^5 +1}{\prod_{i=1}^3x_i} \frac{\partial}{\partial x_4}\right)
\]\\
So we see that {\bf the flow of $V$ is smooth when restricted to the smooth part of $X_{\infty}$.}\\\\
Another way to understand $V$ is to realize that $ds$ is a holomorphic section of $N^*_X$, $(ds)^{-1}$ is a natural holomorphic section of $N_X$. Notice the exact sequence\\
\[
0 \rightarrow T_X \rightarrow T_{\bf P^4}|_X \rightarrow N_X \rightarrow 0
\]\\
With respect to the \k metric $g$ on $T_{\bf P^4}$, the exact sequence has a natural (non-holomorphic) splitting\\
\[
T_{\bf P^4}|_X =T_X \oplus N_X
\]\\
$V$ is just real part of the natural lift of $(ds)^{-1}$ via this splitting. $V$ is singular exactly when $(ds)^{-1}$ is singular or equivalently, when $ds=0$, which corresponds to singular part of $X_{\infty}$. On the other hand, the union of the smooth three-dimensional \l torus fibers of $X_{\infty}$ is exactly the smooth part of $X_{\infty}$. All these 3-torus fibers will be carried to $X_{\psi}$ nicely by the flow of $V$.\\
\subsection{Flow from the singular part of $X_{\infty}$}
Now we will try to understand how the gradient flow of $f$ behaves at singularities of $X_{\infty}$. For this purpose, it is helpful to understand the following example.\\\\
{\bf Example:}
Consider holomorphic function $s(z)=e^{i\theta}z_1z_2\cdots z_n$ on ${\bf C^n}$. $X=\{s=0\}$ is a variety with normal crossing singularities. There is a natural map \\
\[
F:\ {\bf C^n} \rightarrow {\bf R^n_+},\ \ F(z_1,\cdots,z_n) = (|z_1|,\cdots,|z_n|).
\]\\
For $c\in {\bf R^n_+}$, the fiber\\
\[
F^{-1}(c) = \{z:|z_i|=c_i,\ {\rm for}\ i=1,\cdots,n \}
\]\\
is $n$-torus for generic $c$. Let $f=Re(s)$, then with respect to the flat metric on ${\bf C^n}$, we have\\
\[
\nabla f = {\rm Re}\left( s\sum_{i=1}^n \frac{1}{z_i}\frac{\partial}{\partial \bar{z}_i}\right)
\]\\
By our previous argument, $h(z) ={\rm Im}(s)=Im(e^{i\theta}z_1z_2\cdots z_n)$ is invariant under gradient flow of $f$. It is also easy to observe that $\nabla f$ will leave $h_{ij}(z)=|z_i|^2 -|z_j|^2$ invariant. For $c=(c_1,\cdots, c_n)$ such that $\sum c_i =0$, let\\
\[
S_c = \{z: h(z)=0,\ h_{ij} =c_{ij}=c_i-c_j,\ {\rm for}\ i,j=1,\cdots,n \}.
\]\\
It is easy to see that $S_c$ are special Lagrangians in ${\bf C^n}$. Actually this is an example mentioned in Harvey and Lawson's paper. Let $X_r = \{z: s(z)=r\}$, then $L_{c,r} = S_c\cap X_r$ give us a smooth \l $n$-torus fibration of $X_r$ for $r$ real.\\\\
Define\\
\[
\theta_{ij}= \theta_i-\theta_j = \frac{i}{2}\left( \log\frac{\bar{z}_i}{z_i} -\log\frac{\bar{z}_j}{z_j}\right),
\]\\
where $\theta_i$ is the argument of $z_i$. Then\\
\[
\nabla f (\theta_{ij}) =-\frac{1}{2} {\rm Im} \left(\frac{s}{|z_i|^2} - \frac{s}{|z_j|^2}\right)=-\frac{1}{2} \left(\frac{1}{|z_i|^2} - \frac{1}{|z_j|^2}\right)h(z).
\]\\
Recall that the real hypersurface $\{h(z)=0\}$ is invariant under the flow. Above equation implies that when restricted to $\{h(z)=0\}$ $\theta_{ij}$ is invariant. For $c=(c_1,\cdots, c_n)$ such that $\sum c_i =0$, let\\
\[
T_c = \{z: h(z)=0,\ \theta_{ij} =c_{ij}=c_i-c_j,\ {\rm for}\ i,j =1,\cdots,n \}.
\]\\
There is another very illustrative way to write $T_c$.\\
\[
T_c = \{z: \theta_i-\theta_j=c_{ij}=c_i-c_j, \theta+\sum_{i=1}^n\theta_i = 0\}=\{z:\theta_i=c_i-\frac{\theta}{n}\}
\]\\
It is easy to see that $T_c$ are special Lagrangians in ${\bf C^n}$ that is invariant under the flow. $T_{c,r} = T_c\cap X_r$ give us a smooth \l fibration on the horizontal directions of $X_r$ for $r$ real.
\begin{flushright} $\Box$ \end{flushright}
As we know, $X_{\infty}$ has only normal crossing singularities. The above example gave us a rough picture of the local behavior of the gradient flow of $f$ around singularities of $X_{\infty}$ {\bf when denominator of $s$ is non-zero}. For detailed discussion and proof, please refer to \cite{lag2}.\\\\
Finally it remains to analyze the case when denominator of $s$ is zero and $X_{\infty}$ is singular, which is the following set $\Sigma$. Let\\
\[
\Sigma_{ijk} =\{[z]\in {\bf CP^4}| z_i^5 + z_j^5 + z_k^5 =0, z_l=0\ \ {\rm for}\ l\in\{1,2,3,4,5\}\backslash\{i,j,k\}\}.
\]
$\Sigma_{ijk}$ is a genus 6 curve. Let\\
\[
\Sigma = \bigcup_{\{i,j,k\}\in \{1,2,3,4,5\}} \Sigma_{ijk}
\]
We see that $\Sigma = {\rm Sing}(X_{\infty})\cap X_{\psi}$ for any $\psi$.\\\\
In general, points in $X_{\infty}\cap X_{\psi}$ are fixed under the flow. In particular, $\Sigma \subset X_{\infty}\cap X_{\psi}$ is fixed under the flow. Since the vector field $V$ is singular along $\Sigma$, more argument is needed to ensure the flow behave as expected near $\Sigma$. For the argument to work, it is actually necessary to deform the \k form slightly near $\Sigma$. For details, please also refer to \cite{lag2}.\\
\subsection{Lagrangian torus fibration structure}
To understand the \l torus fibration of $X_{\psi}$, it is helpful to first review the torus fibration of $X_{\infty}$. For any subset $I \subset \{1,2,3,4,5\}$, Let\\
\[
D_I = \{z:\ z_i=0,z_j\not= 0,\ {\rm for}\ i\in I,j\in\{1,2,3,4,5\}\backslash I\},
\]\\
 $\Delta_I = F(D_I)\subset \partial\Delta$. Let $|I|$ denote the cardinality of $I$. We have\\
\[
\partial\Delta = \bigcup_{\tiny{\begin{array}{c}I \subset \{1,2,3,4,5\}\\ 0<|I|<5\end{array}}} \Delta_I
\]\\
The fibers over $\Delta_I$ are $4-|I|$ dimensional torus.\\\\
The flow of $V$ moves $X_{\infty}$ to $X_{\psi}$. When $|I|=1$, the flow is diffeomorphism on $D_I$ and moves the corresponding smooth 3-torus fibration in $D_I$ to 3-torus fibration in $X_{\psi}$.\\\\
When $|I|>1$, by above discussion after the example, each point of $D_I\backslash \Sigma$ will be deform to a $|I|-1$ torus under the flow of $V$. Therefore a torus fiber in $D_I$ not intersecting $\Sigma$ will be deformed to a smooth 3-torus \l fiber of $X_{\psi}$ under the flow of $V$.\\\\
Now let's try to understand the singular fibers. Notice that $D_I\cap\Sigma\not= \emptyset$ if and only if $|I|=2,3$. Let $\tilde{\Gamma}_{ijk} = F(\Sigma_{ijk})$. There are 3 types of singular fiber over different portion of\\
\[
\tilde{\Gamma} = \bigcup_{\{i,j,k\}\in \{1,2,3,4,5\}} \tilde{\Gamma}_{ijk}
\]\\
Let $\tilde{\Gamma}^2$ denote the interior of $\tilde{\Gamma}$, and\\
\[
\tilde{\Gamma}^1 = \partial \tilde{\Gamma} \cap \left(\bigcup_{|I|=2}\Delta_I\right),
\]
\[
\tilde{\Gamma}^0 = \partial \tilde{\Gamma} \cap \left(\bigcup_{|I|=3}\Delta_I\right).
\]\\
Then\\
\[
\tilde{\Gamma} = \tilde{\Gamma}^0 \cup \tilde{\Gamma}^1 \cup \tilde{\Gamma}^2.
\]\\
Without loss of generality, let us concentrate on\\
\[
\overline{D}_{45} = \{z: z_4=z_5=0\}
\]\\
with the natural coordinate $x_i=z_i/z_3$ for $i=1,2$. Under this coordinate\\
\[
\Sigma_{123} =\{x=(x_1,x_2)| x_1^5 + x_2^5 + 1 =0\}.
\]\\
$r=(r_1,r_2)=(|x_1|,|x_2|)$ can be thought of as coordinate on $\overline{\Delta}_{45}$. Under this coordinate\\
\[
\tilde{\Gamma}_{123} = \{r=(r_1,r_2)|r_1^5 + r_2^5 \geq 1, r_1^5 +1 \geq r_2^5, r_2^5 +1 \geq r_1^5\}
\]\\
For simplicity, we will omit the index and denote $\Sigma_{123}$ by $\Sigma$ and $\tilde{\Gamma}_{123}$ by $\tilde{\Gamma}$. Then\\
\[
\tilde{\Gamma}^2 = \{r=(r_1,r_2)|r_1^5 + r_2^5 > 1, r_1^5 +1 > r_2^5, r_2^5 +1 > r_1^5 \},
\]
\[
\tilde{\Gamma}^0 = \{(1,0),(0,1)\},
\]
\[
\tilde{\Gamma}^1 = \{r=(r_1,r_2)|r_1^5 + r_2^5 = 1\ {\rm or}\ r_1^5 = r_2^5 + 1\ {\rm or}\ r_2^5 = r_1^5 + 1\}\backslash \tilde{\Gamma}^0.
\]\\
$(x_1,x_2)\rightarrow (e^{\frac{2\pi i}{5}k_1}x_1,e^{\frac{2\pi i}{5}k_2}x_2)$ give us a natural action of ${\bf Z_5\times Z_5}$ on $\Sigma$. Combine with the ${\bf Z_2}$ action $(x_1,x_2)\rightarrow (\bar{x}_1,\bar{x}_2)$ we have\\
\[
\begin{array}{ccc}
{\bf Z_2\times Z_5\times Z_5}&\rightarrow& \Sigma\\
&&\\
&&\downarrow^F\\
&&\\
&&\tilde{\Gamma}\\
\end{array}
\]\\
$F$ is a $50$-sheet covering map over $ \tilde{\Gamma}^2$. ${\bf Z_2}$ is the stablizer of the group action on $F^{-1}(\tilde{\Gamma}^1)\cap\Sigma$. For $p\in \tilde{\Gamma}^1$, $F^{-1}(p)$ is a 2-torus that intersect $\Sigma$ at $25$ points. ${\bf Z_2\times Z_5}$ is the stablizer of the group action on $F^{-1}(\tilde{\Gamma}^0)\cap\Sigma$. For $p\in \tilde{\Gamma}^0$, $F^{-1}(p)$ is a circle that intersect $\Sigma$ at $5$ points.\\\\
From our previous discussion, when $|I|>1$, each point of $D_I\backslash \Sigma$ will be deform to a $|I|-1$ torus under the flow of $V$ and points in $\Sigma$ will not move under the flow of $V$. Let $L_{\infty}\subset D_I \subset X_{\infty}$ be the original fiber of $F$ that move to \l fiber $L\subset X_{\psi}$ under the flow of $V$. Then follow the flow backward we will have a fibration $\pi: L \rightarrow L_{\infty}$. Over $L_{\infty}\backslash \Sigma$, $\pi: L\backslash \Sigma \rightarrow L_{\infty}\backslash \Sigma$ is a $|I|-1$ torus smooth fibration. Over $L_{\infty}\cap \Sigma$, $\pi: L\cap \Sigma \cong L_{\infty}\cap \Sigma$ is an identification.\\\\
It is now easy to see that for $p\in \tilde{\Gamma}^2$, $F^{-1}(p)$ under the flow of $V$ will deform to \l 3-torus with $50$ circles collapsed to $50$ singular points. For $p\in \tilde{\Gamma}^1$, $F^{-1}(p)$ under the flow of $V$ will deform to \l 3-torus with $25$ circles collapsed to $25$ singular points. For $p\in \tilde{\Gamma}^0$, $F^{-1}(p)$ under the flow of $V$ will deform to \l 3-torus with $5$ two-torus collapsed to $5$ singular points. Now we have finished the discussion of all \l fibers of our \l fibration of $X_{\psi}$.\\
\begin{th}
Flow of $V$ will produce a \l fibration $F: X_{\psi} \rightarrow \partial\Delta$. There are 4 types of fibers.\\
(i). For $p\in \partial\Delta\backslash \tilde{\Gamma}$, $F^{-1}(p)$ is a smooth \l 3-torus.\\
(ii). For $p\in \tilde{\Gamma}^2$, $F^{-1}(p)$ is a \l 3-torus with $50$ circles collapsed to $50$ singular points.\\
(iii). For $p\in \tilde{\Gamma}^1$, $F^{-1}(p)$ is a \l 3-torus with $25$ circles collapsed to $25$ singular points.\\
(iv). For $p\in \tilde{\Gamma}^0$, $F^{-1}(p)$ is a \l 3-torus with $5$ two-torus collapsed to $5$ singular points.
\end{th}
\begin{flushright} $\Box$ \end{flushright}

\se{Monodromy of expected \sl fibration}
\subsection{Introduction and assumptions}
In this section we will explore what the SYZ \sl fibration of $X_{\psi}$ should look like if it exists. Recall from Section 2 that $X_{\infty}$ is a union of five ${\bf P^3}$. There is a natural degenerate  $T^3$ fibration structure for $X_{\infty}$. Let $\{P_i\}_{i=1}^5$ be five points in ${\bf R^4}$ that are in general position. Consider the natural map $F: \ {\bf P^4}\longrightarrow {\bf R^4}$.\\
\[
F([z])=\sum_{k=1}^5 \frac{|z_k|^2}{\sum_{i=1}^5 |z_i|^2}P_k
\]\\
$\Delta={\rm Image}(F)$ is a 4-simplex. $X_{\infty}$ is naturally fibered over $\partial \Delta$ via this map $F$ with general fiber being $T^3$. This is precisely the $T^3$ special Lagrangian fibration for $X_{\infty}$ as indicated by SYZ construction.\\\\
When $\psi$ is large, we expect that $X_{\psi}$ will also possess a special Largrangian $T^3$ fibration with base identified with $\partial \Delta$ which topologically is an $S^3$. There are still serious analysis and geometric works to be done to totally justify the special Lagrangian fibration. We will discuss in this section that {\bf IF} such special Lagrangian fibration exist on $X_{\psi}$, what should be its expected topological and geometrical structures.\\\\
Our discussion is based on two principle assumptions:\\
(1) \sl fibration on $X_{\psi}$ is a deformation of \sl fibration on $X_{\infty}$. \\
(2) Singular locus of the \sl fibration on $X_{\psi}$ is of codimension 2. \\\\
(1) is very natural, because we expect the \cy metric on $X_{\psi}$ suitably rescaled (keeping fiber class constant volume) will approach the standard \cy metric on $X_{\infty}$ when $\psi$ approach infinity. (2) is also very reasonable given the structure of elliptic fibration of complex surface that correspond to dimension $n=2$ situation.\\\\
For $\psi$ large, $X_{\psi}$ will approach $X_{\infty}=\cup_{k=1}^5 D_k$, where $D_k=\{z_k=0\}$. Let us consider the part of the $X_{\psi}$ that is close to $D_5$. (We denote it as $U_5$, for example we may take $U_5$ as inverse image in $X_{\psi}$ of an open set in $D_5$ that stay away from $D_k$ for $k\not=5$ by the following map $\pi_5$.) Then the restriction of the projection\\
\[
\pi_5:\ U_5 \longrightarrow D_5,\ \ \pi_5([z_1,z_2,z_3,z_4,z_5])=[z_1,z_2,z_3,z_4,0]
\]
will identify $U_5$ with an open set in $D_5$, which at the same time will carry over the $T^3$ fibration to $U_5$ from $D_5$. When $\psi$ is large, this fibration should be very close to the special Lagrangian $T^3$ fibration for $U_5$. When we discuss the topology aspect, it will be sufficient to use the induced fibration instead of the special Lagrangian fibration when staying away from intersections of $D_k$'s, where the special Lagrangian fibration of $X_{\infty}$ degenerate. We will use these identifications to compute monodromy of the special Lagrangian fibration in the following.\\
\subsection{One forms and one cycles}
We first fix some notation and introduce some construction. On $D_5$, let $\gamma^1_{52}$ be the circles determined by $\{z_5=0,\ |\frac{z_1}{z_2}| = c_1,\ \frac{z_3}{z_2}=c_2,\ \frac{z_4}{z_2}=c_3\}$. We will also use it to denote circles carried over to $U_5$. In general, we have $\gamma_{ik}^j$ for $\{i,j,k\} \subset \{1,2,3,4,5\}$. Understanding monodromy is equivalent to understanding transformations among $\{\gamma_{ik}^j\}$. It is easy to check that\\
\begin{eqnarray}
\gamma_{ik}^j &=& \gamma_{il}^j,\ \ \ \ \ {\rm for}\ \{i,j,k,l\}\subset \{1,2,3,4,5\}.\label{ba}\\
\gamma_{ik}^j &=& -\sum_{l(\not=i,j)=1}^5\gamma_{ij}^l,\ \ {\rm for}\ \{i,j,k\}\subset \{1,2,3,4,5\}.\nonumber
\end{eqnarray}
On ${\bf P^4}$ we can introduce meromorphic 1-forms\\
\[
\alpha_{ij}= d(\log \frac{z_i}{z_j}), \ \ i \not=j.
\]
They have the simple relations\\
\[
\alpha_{ij}= \alpha_{ik}+ \alpha_{kj}\ {\rm for}\ \{i,j,k\}\subset \{1,2,3,4,5\}.
\]
We also use the same notation for the restriction to $X_{\psi}$. $\alpha_{ij}$ have pole along $D_j\cap X_{\psi}$ on $X_{\psi}$, where\\
\[
D_j\cap X_{\psi}=\left\{z\left| z_j=0,\ \sum_{i(\not=j)=1}^5 z_i^5 = 0\right. \right\}
\]
For this reason we introduce\\
\[
U_j^i = \left\{z\in U_j \left|\ \sum_{k(\not= i,j)=1}^5\left|\frac{z_k}{z_i}\right|^5<1\right. \right\}.
\]
On $U_j^i$, all $\alpha_{kl}$ are regular. On $U_i^j$ we have circles \\
\[
\gamma_{ij}^k, \ {\rm for}\ k\in \{1,2,3,4,5\}\backslash \{i,j\}
\]
and 1-forms\\
\[
\alpha_{kj},\ {\rm for}\ k\in \{1,2,3,4,5\}\backslash \{j\}
\]
It is easy to check that\\
\[
<\gamma_{ij}^k, \alpha_{lj}> = \delta_{kl}\ {\rm for}\ \{i,j,k,l\}\subset \{1,2,3,4,5\}
\]
and\\
\[
<\gamma_{ij}^k, \alpha_{ij}> = -1\ {\rm for}\ \{i,j,k\}\subset \{1,2,3,4,5\}
\]
Use these relations for $U_i^j$ and $U_k^j$, we will get\\
\begin{eqnarray}
\gamma_{ij}^k&=&-\gamma_{kj}^i,\ \ \ \ {\rm for}\ \{i,j,k\}\subset \{1,2,3,4,5\}\label{bb}\\
\gamma_{ij}^l&=&-\gamma_{kj}^i+\gamma_{kj}^l,\ \ {\rm for}\ \{i,j,k,l\}\subset \{1,2,3,4,5\}\nonumber
\end{eqnarray}
From (\ref{ba}) we can see that without confusion, we may denote $\gamma_{ik}^j$ by $\gamma_i^j$. then the relations (\ref{ba}) and (\ref{bb}) can be rewrite together as\\
\begin{eqnarray}
\sum_{j(\not=i)=1}^5\gamma_i^j&=&0,\ \ {\rm for}\ \{i\}\subset \{1,2,3,4,5\}\nonumber\\
\gamma_i^j+\gamma_j^i&=&0,\ \ \ {\rm for}\ \{i,j\}\subset \{1,2,3,4,5\}\label{bc}\\
\gamma_i^j+\gamma_j^k+\gamma_k^i&=&0,\ \ \ {\rm for}\ \{i,j,k\}\subset \{1,2,3,4,5\}\nonumber
\end{eqnarray}
Let\\
\[
U^i = \left\{z\in X_{\psi} \left|\ \sum_{k(\not= i)=1}^5\left|\frac{z_k}{z_i}\right|^5<1\right. \right\}.
\]
Then $\alpha_{ki}$, for $k\in \{1,2,3,4,5\}\backslash \{i\}$ are regular on $U^i$. So there are no monodromy among $U^i_j\subset U^i$ for $j\in \{1,2,3,4,5\}\backslash \{i\}$.\\\\
On the other hand, in $U_j$ the $T^3$ fibration is regular. So there are also no monodromy among $U^i_j\subset U_j$ for $j\in \{1,2,3,4,5\}\backslash \{i\}$. From these discussions we can see that the discriminant locus of the fibration (where the $T^3$ fibration is singular) is topologically a graph $\Gamma$ in $\partial \Delta$. Vertices of $\Gamma$ are $P_{ij}$ (baricenter for 2-simplex $\Delta_{ij}$) for $\{i,j\}\subset \{1,2,3,4,5\}$ and $P_{ijk}$ (baricenter for 3-simplex $\Delta_{ijk}$) for $\{i,j,k\}\subset \{1,2,3,4,5\}$. Legs of $\Gamma$ are $\Gamma_{ij}^k$ which connects $P_{ijk}$ and $P_{ij}$ for $\{i,j,k\}\subset \{1,2,3,4,5\}$.\\\\
\begin{figure}[h]
\begin{center}
\leavevmode
\hbox{%
\epsfxsize=4in
\epsffile{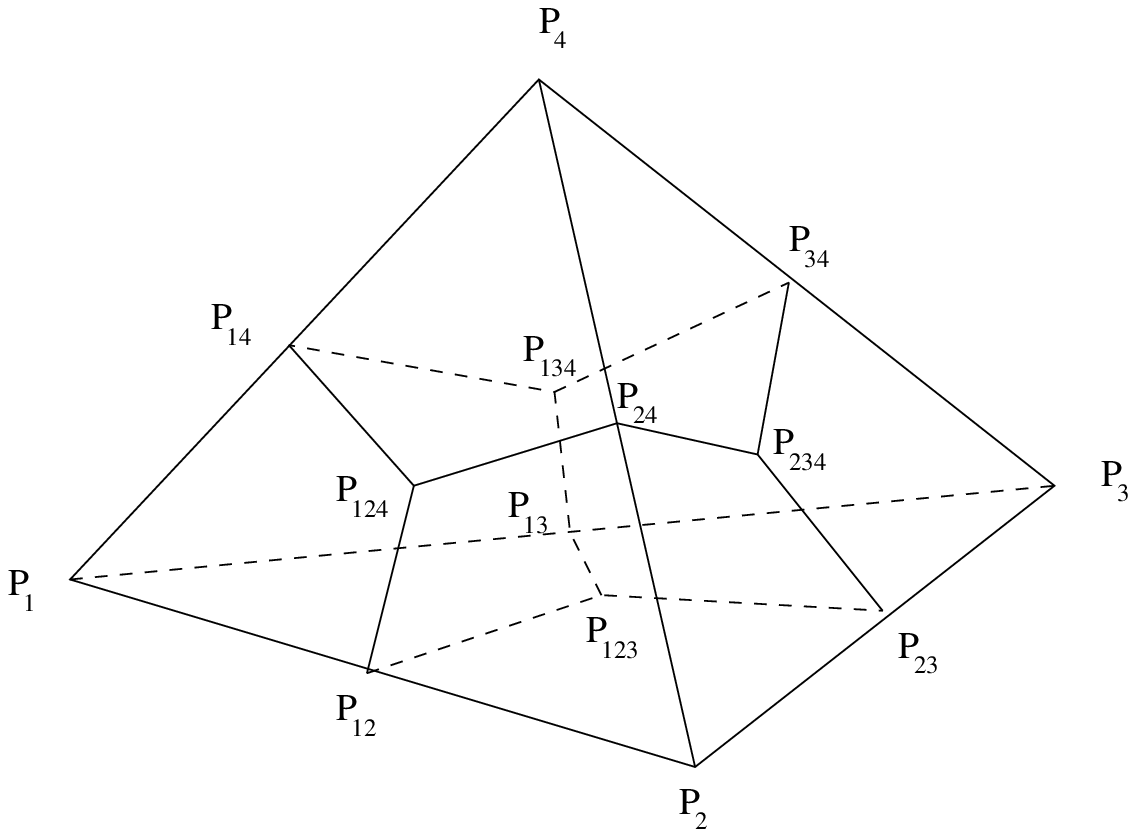}}
\end{center}
\caption{}
\end{figure}\\
\subsection{Monodromy around the legs of $\Gamma$}
We would like to compute monodromy $T_{ij}^k$ around $\Gamma_{ij}^k$. To simplify the notation without loss of generality, we will compute $T_{24}^3$. This amounts to compute the transformation:\\
\[
U_5^4\rightarrow U_5^2\rightarrow U_1^2\rightarrow U_1^4\rightarrow U_5^4
\]
Let $S_i^{jk}$ denote the transformation: $U_i^j\rightarrow U_i^k$ and $S^i_{jk}$ denote the transformation: $U^i_j\rightarrow U^i_k$ for $\{i,j,k\}\subset \{1,2,3,4,5\}$. Then\\
\[
T_{24}^3= S_{15}^4 S_1^{24} S_{51}^2 S_5^{42}
\]
For $U_5^4\rightarrow U_5^2$ we have\\
\[
S_5^{42}\left(
\begin{array}{ccc}
\gamma_{54}^1 & \gamma_{54}^2 & \gamma_{54}^3
\end{array}
\right) = \left(
\begin{array}{ccc}
\gamma_{52}^1 & \gamma_{52}^3 & \gamma_{52}^4
\end{array}
\right) \left(
\begin{array}{ccc}
1 & -1 & 0\\
0 & -1 & 1\\
0 & -1 & 0
\end{array}
\right)
\]
For $U_5^2\rightarrow U_1^2$ we have\\
\[
S_{51}^2\left(
\begin{array}{ccc}
\gamma_{52}^1 & \gamma_{52}^3 & \gamma_{52}^4
\end{array}
\right) = \left(
\begin{array}{ccc}
\gamma_{12}^3 & \gamma_{12}^4 & \gamma_{12}^5
\end{array}
\right) \left(
\begin{array}{ccc}
0 & 1 & 0\\
0 & 0 & 1\\
-1& -1& -1
\end{array}
\right)
\]
For $U_1^2\rightarrow U_1^4$ we have\\
\[
S_1^{24}\left(
\begin{array}{ccc}
\gamma_{12}^3 & \gamma_{12}^4 & \gamma_{12}^5
\end{array}
\right) = \left(
\begin{array}{ccc}
\gamma_{14}^2 & \gamma_{14}^3 & \gamma_{14}^5
\end{array}
\right) \left(
\begin{array}{ccc}
0 & -1 & 0\\
1 & -1 & 0\\
0 & -1 & 1
\end{array}
\right)
\]
For $U_1^4\rightarrow U_5^4$ we have\\
\[
S_{15}^2\left(
\begin{array}{ccc}
\gamma_{14}^2 & \gamma_{14}^3 & \gamma_{14}^5
\end{array}
\right) = \left(
\begin{array}{ccc}
\gamma_{54}^1 & \gamma_{54}^2 & \gamma_{54}^3
\end{array}
\right) \left(
\begin{array}{ccc}
-1&-1 &-1\\
1 & 0 & 0\\
0 & 1 & 0
\end{array}
\right)
\]
Therefore\\
\[
T_{24}^3= S_{15}^4 S_1^{24} S_{51}^2 S_5^{42}
\]
\[
=\left(
\begin{array}{ccc}
-1&-1 &-1\\
1 & 0 & 0\\
0 & 1 & 0
\end{array}
\right)
\left(
\begin{array}{ccc}
0 & -1 & 0\\
1 & -1 & 0\\
0 & -1 & 1
\end{array}
\right)
\left(
\begin{array}{ccc}
0 & 1 & 0\\
0 & 0 & 1\\
-1& -1& -1
\end{array}
\right)
\left(
\begin{array}{ccc}
1 & -1 & 0\\
0 & -1 & 1\\
0 & -1 & 0
\end{array}
\right)
\]
\[
=\left(
\begin{array}{ccc}
1 & -5 & 0\\
0 &  1 & 0\\
0 &  0 & 1
\end{array}
\right)
\]
and\\
\[
T_{24}^3 \left(
\begin{array}{ccc}
\gamma_{54}^1 & \gamma_{54}^2 & \gamma_{54}^3
\end{array}
\right)=\left(
\begin{array}{ccc}
\gamma_{54}^1 & \gamma_{54}^2 & \gamma_{54}^3
\end{array}
\right)\left(
\begin{array}{ccc}
1 & -5 & 0\\
0 &  1 & 0\\
0 &  0 & 1
\end{array}
\right)
\]
In general for $\{i,j,k,l,m\}=\{1,2,3,4,5\}$ representing the same orientation with the written order\\
\[
T_{jl}^k \left(
\begin{array}{ccc}
\gamma_{m}^i & \gamma_{m}^j & \gamma_{m}^k
\end{array}
\right)=\left(
\begin{array}{ccc}
\gamma_{m}^i & \gamma_{m}^j & \gamma_{m}^k
\end{array}
\right)\left(
\begin{array}{ccc}
1 & -5 & 0\\
0 &  1 & 0\\
0 &  0 & 1
\end{array}
\right),
\]\\
where the monodromy is computed along the path\\
\[
U_m^l\rightarrow U_m^j\rightarrow U_i^j\rightarrow U_i^l\rightarrow U_m^l.
\]\\
If orientation is different, then $-5$ should be replaced by $5$.\\
\subsection{Monodromy around baricenter of 2-cell}
Now we compute the monodromy around a baricenter of a 2-cell.\\\\
Let $S_m^{kl}$ denote the transformation: $U_m^k\rightarrow U_m^l$\\
\[
S_m^{kl}\left(
\begin{array}{ccc}
\gamma_m^i & \gamma_m^j & \gamma_m^l
\end{array}
\right) = \left(
\begin{array}{ccc}
\gamma_m^i & \gamma_m^j & \gamma_m^k
\end{array}
\right) \left(
\begin{array}{ccc}
1 & 0 & -1\\
0 & 1 & -1\\
0 & 0 & -1
\end{array}
\right)
\]\\
and $S^k_{ml}$ denote the transformation: $U^k_m\rightarrow U^k_l$\\
\[
S_{ml}^k\left(
\begin{array}{ccc}
\gamma_m^i & \gamma_m^j & \gamma_m^l
\end{array}
\right) = \left(
\begin{array}{ccc}
\gamma_l^i & \gamma_l^j & \gamma_l^m
\end{array}
\right) \left(
\begin{array}{ccc}
1 & 0 & 0\\
0 & 1 & 0\\
-1& -1& -1
\end{array}
\right)
\]\\
For $(i,j,k,l,m)=(1,4,2,3,5)$\\
\[
T_{43}^2 \left(
\begin{array}{ccc}
\gamma_5^1 & \gamma_5^4 & \gamma_5^2
\end{array}
\right)=\left(
\begin{array}{ccc}
\gamma_5^1 & \gamma_5^4 & \gamma_5^2
\end{array}
\right)\left(
\begin{array}{ccc}
1 & -5 & 0\\
0 &  1 & 0\\
0 &  0 & 1
\end{array}
\right)
\]\\
For $(i,j,k,l,m)=(1,3,4,2,5)$\\
\[
T_{32}^4 \left(
\begin{array}{ccc}
\gamma_5^1 & \gamma_5^3 & \gamma_5^4
\end{array}
\right)=\left(
\begin{array}{ccc}
\gamma_5^1 & \gamma_5^3 & \gamma_5^4
\end{array}
\right)\left(
\begin{array}{ccc}
1 & -5 & 0\\
0 &  1 & 0\\
0 &  0 & 1
\end{array}
\right)
\]\\
Write under the same basis\\
\[
T_{43}^2 \left(
\begin{array}{ccc}
\gamma_{5}^1 & \gamma_{5}^2 & \gamma_{5}^3
\end{array}
\right)=
T_{43}^2 \left(
\begin{array}{ccc}
\gamma_{5}^1 & \gamma_{5}^2 & \gamma_{5}^4
\end{array}
\right)\left(
\begin{array}{ccc}
1 & 0 & -1\\
0 & 1 & -1\\
0 & 0 & -1
\end{array}
\right)
\]\\
\[
=\left(
\begin{array}{ccc}
\gamma_{5}^1 & \gamma_{5}^2 & \gamma_{5}^4
\end{array}
\right)\left(
\begin{array}{ccc}
1 & 0 & -5\\
0 & 1 & 0\\
0 & 0 & 1
\end{array}
\right)\left(
\begin{array}{ccc}
1 & 0 & -1\\
0 & 1 & -1\\
0 & 0 & -1
\end{array}
\right)
\]\\
\[
=\left(
\begin{array}{ccc}
\gamma_{5}^1 & \gamma_{5}^2 & \gamma_{5}^3
\end{array}
\right)\left(
\begin{array}{ccc}
1 & 0 & -1\\
0 & 1 & -1\\
0 & 0 & -1
\end{array}
\right)\left(
\begin{array}{ccc}
1 & 0 & -5\\
0 & 1 & 0\\
0 & 0 & 1
\end{array}
\right)\left(
\begin{array}{ccc}
1 & 0 & -1\\
0 & 1 & -1\\
0 & 0 & -1
\end{array}
\right)
\]\\
which gives\\
\[
T_{43}^2 \left(
\begin{array}{ccc}
\gamma_{5}^1 & \gamma_{5}^2 & \gamma_{5}^3
\end{array}
\right)=\left(
\begin{array}{ccc}
\gamma_{5}^1 & \gamma_{5}^2 & \gamma_{5}^3
\end{array}
\right)\left(
\begin{array}{ccc}
1 &  0 & 5\\
0 &  1 & 0\\
0 &  0 & 1
\end{array}
\right)
\]\\
\[
T_{32}^4 \left(
\begin{array}{ccc}
\gamma_{5}^1 & \gamma_{5}^2 & \gamma_{5}^3
\end{array}
\right)=
T_{32}^4 \left(
\begin{array}{ccc}
\gamma_{5}^1 & \gamma_{5}^3 & \gamma_{5}^4
\end{array}
\right)\left(
\begin{array}{ccc}
1 & -1 & 0\\
0 & -1 & 1\\
0 & -1 & 0
\end{array}
\right)
\]\\
\[
=\left(
\begin{array}{ccc}
\gamma_{5}^1 & \gamma_{5}^3 & \gamma_{5}^4
\end{array}
\right)\left(
\begin{array}{ccc}
1 & -5 & 0\\
0 &  1 & 0\\
0 &  0 & 1
\end{array}
\right)\left(
\begin{array}{ccc}
1 & -1 & 0\\
0 & -1 & 1\\
0 & -1 & 0
\end{array}
\right)
\]\\
\[
=\left(
\begin{array}{ccc}
\gamma_{5}^1 & \gamma_{5}^2 & \gamma_{5}^3
\end{array}
\right)\left(
\begin{array}{ccc}
1 & 0 & -1\\
0 & 0 & -1\\
0 & 1 & -1
\end{array}
\right)\left(
\begin{array}{ccc}
1 & -5 & 0\\
0 &  1 & 0\\
0 &  0 & 1
\end{array}
\right)\left(
\begin{array}{ccc}
1 & -1 & 0\\
0 & -1 & 1\\
0 & -1 & 0
\end{array}
\right)
\]\\
which gives\\
\[
T_{32}^4 \left(
\begin{array}{ccc}
\gamma_{5}^1 & \gamma_{5}^2 & \gamma_{5}^3
\end{array}
\right)=\left(
\begin{array}{ccc}
\gamma_{5}^1 & \gamma_{5}^2 & \gamma_{5}^3
\end{array}
\right)\left(
\begin{array}{ccc}
1 & 5 & -5\\
0 & 1 &  0\\
0 & 0 &  1
\end{array}
\right)
\]\\
From these expressions, it is clear that $T_{32}^4$, $T_{24}^3$ and $T_{43}^2$ commute with each other and\\
\[
T_{32}^4 T_{24}^3 T_{43}^2= Id.
\]\\
They determine a natural filtration:\\
\[
{\cal W}^0 \subset {\cal W}^1 = H_1(T^3, {\bf Z})
\]\\
with ${\cal W}^0$ generated by the vanishing circle $\gamma_{5}^1$, which is the common vanishing cycle of $T_{32}^4$, $T_{24}^3$ and $T_{43}^2$. Recall that the $T^3$ fibration of $X_{\infty}$ over $\partial \Delta$ degenerate to be a $T^2$ fibration over $\Delta _{15}$. An interesting fact is that this $T^2$ is exactly the quotient of $T^3$ by $\gamma_{5}^1$.\\
\subsection{Monodromy around baricenter of 1-cell}
We can similarly compute the monodromies around a midpoint of a 1-edge, for instance, $P_{24}$. We have for $(i,j,k,l,m)=(5,2,1,4,3)$\\
\[
T_{24}^1 \left(
\begin{array}{ccc}
\gamma_{3}^5 & \gamma_{3}^2 & \gamma_{3}^1
\end{array}
\right)=\left(
\begin{array}{ccc}
\gamma_{3}^5 & \gamma_{3}^2 & \gamma_{3}^1
\end{array}
\right)\left(
\begin{array}{ccc}
1 & -5 & 0\\
0 &  1 & 0\\
0 &  0 & 1
\end{array}
\right)
\]\\
and for $(i,j,k,l,m)=(3,2,5,4,1)$\\
\[
T_{24}^5 \left(
\begin{array}{ccc}
\gamma_{1}^3 & \gamma_{1}^2 & \gamma_{1}^5
\end{array}
\right)=\left(
\begin{array}{ccc}
\gamma_{1}^3 & \gamma_{1}^2 & \gamma_{1}^5
\end{array}
\right)\left(
\begin{array}{ccc}
1 & -5 & 0\\
0 &  1 & 0\\
0 &  0 & 1
\end{array}
\right)
\]\\
write under the same basis $(\gamma_{5}^1, \gamma_{5}^2, \gamma_{5}^3)$.\\
\[
T_{24}^1 \left(
\begin{array}{ccc}
\gamma_{5}^1 & \gamma_{5}^2 & \gamma_{5}^3
\end{array}
\right)=
T_{24}^1 \left(
\begin{array}{ccc}
\gamma_{3}^5 & \gamma_{3}^2 & \gamma_{3}^1
\end{array}
\right)\left(
\begin{array}{ccc}
-1 & -1 & -1\\
 0 &  1 &  0\\
 1 &  0 &  0
\end{array}
\right)
\]\\
\[
=\left(
\begin{array}{ccc}
\gamma_{3}^5 & \gamma_{3}^2 & \gamma_{3}^1
\end{array}
\right)\left(
\begin{array}{ccc}
1 & -5 & 0\\
0 &  1 & 0\\
0 &  0 & 1
\end{array}
\right)\left(
\begin{array}{ccc}
-1&-1 & -1\\
0 & 1 &  0\\
1 & 0 &  0
\end{array}
\right)
\]\\
\[
=\left(
\begin{array}{ccc}
\gamma_{5}^1 & \gamma_{5}^2 & \gamma_{5}^3
\end{array}
\right)\left(
\begin{array}{ccc}
0 &  0 & 1\\
0 &  1 & 0\\
-1& -1 &-1
\end{array}
\right)\left(
\begin{array}{ccc}
1 & -5 & 0\\
0 &  1 & 0\\
0 &  0 & 1
\end{array}
\right)\left(
\begin{array}{ccc}
-1&-1 & -1\\
0 & 1 &  0\\
1 & 0 &  0
\end{array}
\right)
\]\\
which gives\\
\[
T_{24}^1 \left(
\begin{array}{ccc}
\gamma_{5}^1 & \gamma_{5}^2 & \gamma_{5}^3
\end{array}
\right)=\left(
\begin{array}{ccc}
\gamma_{5}^1 & \gamma_{5}^2 & \gamma_{5}^3
\end{array}
\right)\left(
\begin{array}{ccc}
1 &  0 & 0\\
0 &  1 & 0\\
0 &  5 & 1
\end{array}
\right)
\]\\
\[
T_{24}^5 \left(
\begin{array}{ccc}
\gamma_{5}^1 & \gamma_{5}^2 & \gamma_{5}^3
\end{array}
\right)=
T_{24}^5 \left(
\begin{array}{ccc}
\gamma_{1}^3 & \gamma_{14}^2 & \gamma_{1}^5
\end{array}
\right)\left(
\begin{array}{ccc}
 0 &  0 &  1\\
 0 &  1 &  0\\
-1 & -1 & -1
\end{array}
\right)
\]\\
\[
=\left(
\begin{array}{ccc}
\gamma_{1}^3 & \gamma_{1}^2 & \gamma_{1}^5
\end{array}
\right)\left(
\begin{array}{ccc}
1 & -5 & 0\\
0 &  1 & 0\\
0 &  0 & 1
\end{array}
\right)\left(
\begin{array}{ccc}
 0 &  0 &  1\\
 0 &  1 &  0\\
-1 & -1 & -1
\end{array}
\right)
\]\\
\[
=\left(
\begin{array}{ccc}
\gamma_{5}^1 & \gamma_{5}^2 & \gamma_{5}^3
\end{array}
\right)\left(
\begin{array}{ccc}
-1& -1 &-1\\
 0&  1 & 0\\
 1&  0 & 0
\end{array}
\right)\left(
\begin{array}{ccc}
1 & -5 & 0\\
0 &  1 & 0\\
0 &  0 & 1
\end{array}
\right)\left(
\begin{array}{ccc}
 0 &  0 &  1\\
 0 &  1 &  0\\
-1 & -1 & -1
\end{array}
\right)
\]\\
which gives\\
\[
T_{24}^5 \left(
\begin{array}{ccc}
\gamma_{5}^1 & \gamma_{5}^2 & \gamma_{5}^3
\end{array}
\right)=\left(
\begin{array}{ccc}
\gamma_{5}^1 & \gamma_{5}^2 & \gamma_{5}^3
\end{array}
\right)\left(
\begin{array}{ccc}
1 &  5 & 0\\
0 &  1 & 0\\
0 & -5 & 1
\end{array}
\right)
\]\\
From these expressions, it is clear that $T_{24}^1$, $T_{24}^3$ and $T_{24}^5$ commute with each other and\\
\[
T_{24}^1 T_{24}^3 T_{24}^5= Id.
\]\\
They determine a natural filtration:\\
\[
{\cal W}^0 \subset {\cal W}^1 = H_1(T^3, {\bf Z})
\]\\
with ${\cal W}^0$ generated by the vanishing circle $\gamma_{5}^1$ and $\gamma_{5}^3$. Recall that the $T^3$ fibration of $X_{\infty}$ over $\partial \Delta$ degenerate to be a $T^1$ fibration over $\Delta _{135}$. An interesting fact is that this $T^1$ is exactly the quotient of $T^3$ by $\gamma_{5}^1\times \gamma_{5}^3$.\\\\

\se{Geometry of the singular fibers}
With monodromy information in mind, we would like to discuss possible structure of singular special Lagrangian fibers. We will start by describing a possible model for the structure of \sl fibration, especially the singular fibers, that is of conjectural nature. Then we will use our construction of \l fibration to give some partial justification.\\\\
Suppose that we have a \sl fibration $F: X_{\psi} \rightarrow \partial \Delta$. Then the monodromy information will suggest that it is likely that $F$ is smooth fibration over $\partial \Delta\backslash \Gamma$, where the singular locus $\Gamma$ is topologically a graph in $\partial \Delta$. Vertices of $\Gamma$ are $P_{ij}$ (baricenter for 2-simplex $\Delta_{ij}$) for $\{i,j\}\subset \{1,2,3,4,5\}$ and $P_{ijk}$ (baricenter for 3-simplex $\Delta_{ijk}$) for $\{i,j,k\}\subset \{1,2,3,4,5\}$. Legs of $\Gamma$ are $\Gamma_{ij}^k$ which connects $P_{ijk}$ and $P_{ij}$ for $\{i,j,k\}\subset \{1,2,3,4,5\}$.\\\\
In general for $\{i,j,k,l,m\}=\{1,2,3,4,5\}$ representing the same orientation with the written order, the monodromy around $\Gamma_{ij}^k$ computed along the path\\
\[
U_m^j\rightarrow U_m^i\rightarrow U_l^i\rightarrow U_l^j\rightarrow U_m^j.
\]\\
can be written as\\
\[
T_{ij}^k \left(
\begin{array}{ccc}
\gamma_{m}^l & \gamma_{m}^i & \gamma_{m}^k
\end{array}
\right)=\left(
\begin{array}{ccc}
\gamma_{m}^l & \gamma_{m}^i & \gamma_{m}^k
\end{array}
\right)\left(
\begin{array}{ccc}
1 & -5 & 0\\
0 &  1 & 0\\
0 &  0 & 1
\end{array}
\right).
\]\\
This suggests that singular fibers on $\Gamma_{ij}^k$ should be of $I_5 \times S^1$ type, where $I_5$ is a Kodaira type elliptic singular fiber as degeneration of $\gamma_{m}^l \times \gamma_{m}^i$ with vanishing cycle $\gamma_{m}^l$, and the $S^1$ corresponds to $\gamma_{m}^k$. This is a general type of singular fibers. Let's call it type $I_5$ singular fiber. Since Type $I_5$ singular fibers have $S^1$ factor, they give no contribution to Euler number.\\\\
On the other hand, around $P_{ijk}$, monodromies are\\
\[
T_{ij}^k \left(
\begin{array}{ccc}
\gamma_m^l & \gamma_m^i & \gamma_m^k
\end{array}
\right)=\left(
\begin{array}{ccc}
\gamma_m^l & \gamma_m^i & \gamma_m^k
\end{array}
\right)\left(
\begin{array}{ccc}
1 & -5 & 0\\
0 &  1 & 0\\
0 &  0 & 1
\end{array}
\right).
\]
\[
T_{jk}^i \left(
\begin{array}{ccc}
\gamma_m^l & \gamma_m^i & \gamma_m^k
\end{array}
\right)=\left(
\begin{array}{ccc}
\gamma_m^l & \gamma_m^i & \gamma_m^k
\end{array}
\right)\left(
\begin{array}{ccc}
1 &  0 & 5\\
0 &  1 & 0\\
0 &  0 & 1
\end{array}
\right).
\]
\[
T_{ki}^j \left(
\begin{array}{ccc}
\gamma_m^l & \gamma_m^i & \gamma_m^k
\end{array}
\right)=\left(
\begin{array}{ccc}
\gamma_m^l & \gamma_m^i & \gamma_m^k
\end{array}
\right)\left(
\begin{array}{ccc}
1 & 5 &-5\\
0 & 1 & 0\\
0 & 0 & 1
\end{array}
\right).
\]
They determine the natural filtration:\\
\[
{\cal W}^0 \subset {\cal W}^1 = H_1(T^3, {\bf Z})
\]
with ${\cal W}^0$ generated by the common vanishing cycle $\gamma_m^l$. By the symmetric property, the singular fiber over $P_{ijk}$ can be described as follows. Let $\Gamma_{ijk}$ be the graph in $T^2\cong T^3/\gamma_{m}^l$ as discribed by the following picture, where $\{lm\}=\overline{\{ijk\}}$. Then the singular fiber over $P_{ijk}$ can be identified topologically as $T^3$ collapsing circles $(\gamma_{m}^l)$ over $\Gamma_{ijk}$ to points. We call this kind of fibers type $II_{5\times 5}$. It is easy to see that the Euler number of a type $II_{5\times 5}$ fiber is $-25$.\\
\begin{figure}[h]
\begin{center}
\leavevmode
\hbox{%
\epsfxsize=3in
\epsffile{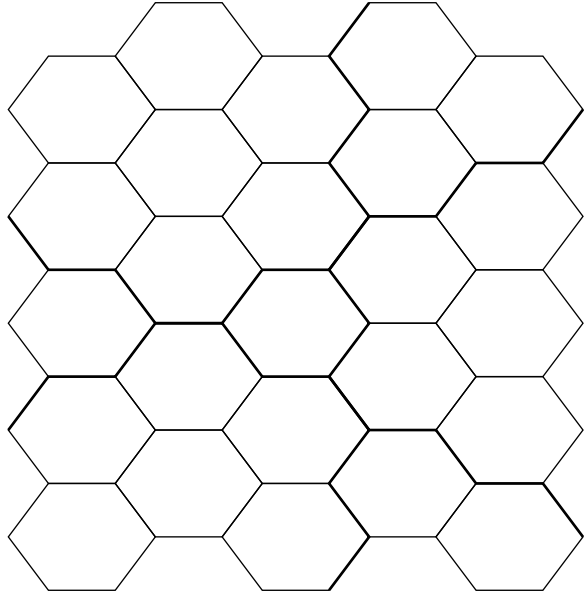}}
\end{center}
\caption{Type $II_{5\times 5}$ fiber}
\end{figure}\\
Finally we are left to determine the singular fiber over $P_{ij}$. Around $P_{ij}$ the monodromies are\\
\[
T_{ij}^k \left(
\begin{array}{ccc}
\gamma_m^l & \gamma_m^i & \gamma_m^k
\end{array}
\right)=\left(
\begin{array}{ccc}
\gamma_m^l & \gamma_m^i & \gamma_m^k
\end{array}
\right)\left(
\begin{array}{ccc}
1 & -5 & 0\\
0 &  1 & 0\\
0 &  0 & 1
\end{array}
\right).
\]
\[
T_{ij}^l \left(
\begin{array}{ccc}
\gamma_m^l & \gamma_m^i & \gamma_m^k
\end{array}
\right)=\left(
\begin{array}{ccc}
\gamma_m^l & \gamma_m^i & \gamma_m^k
\end{array}
\right)\left(
\begin{array}{ccc}
1 & 0 & 0\\
0 & 1 & 0\\
0 & 5 & 1
\end{array}
\right).
\]
\[
T_{ij}^l \left(
\begin{array}{ccc}
\gamma_m^l & \gamma_m^i & \gamma_m^k
\end{array}
\right)=\left(
\begin{array}{ccc}
\gamma_m^l & \gamma_m^i & \gamma_m^k
\end{array}
\right)\left(
\begin{array}{ccc}
1 &  5 & 0\\
0 &  1 & 0\\
0 & -5 & 1
\end{array}
\right).
\]
They determine the natural filtration:\\
\[
{\cal W}^0 \subset {\cal W}^1 = H_1(T^3, {\bf Z})
\]
with ${\cal W}^0$ generated by the vanishing circles $\gamma_m^l$ and $\gamma_m^k$. The singular fiber over $P_{ij}$ is sort of generalization of $I_5$ singularity to $3$-dimension. Topologically equivalent to $T^3$ collapsing five of sub-$T^2$ equivalent to $\gamma_m^l\times \gamma_m^k$. It looks like five $3$-dimensional pseudomanifolds linked into a ``necklace'' via their singular points. Each $3$-dimensional pseudomanifold is a $T^2\times [0,1]$ collapsing the two boundaries to two singular points. An illustration of this singular fiber is in the following picture. We will call these type $III_5$ singular fibers. A type $III_5$ singular fiber has Euler number $5$.\\
\begin{figure}[h]
\begin{center}
\leavevmode
\hbox{%
\epsfxsize=4in
\epsffile{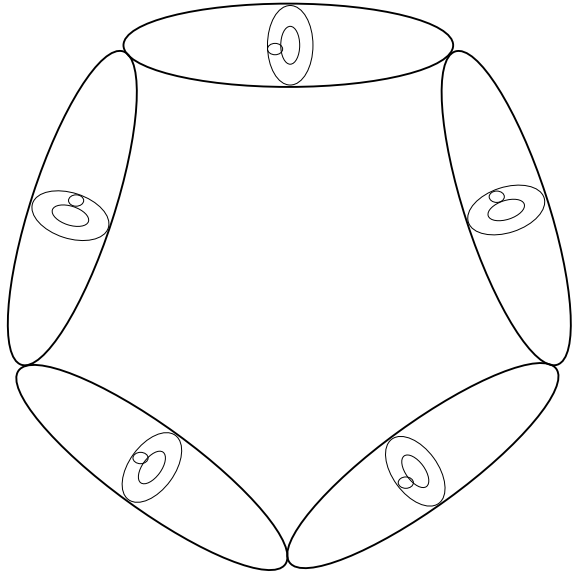}}
\end{center}
\caption{Type $III_5$ fiber}
\end{figure}\\
All together we have $10$ type $II_{5\times 5}$ and $10$ type $III_5$ singular fibers, which give total Euler number\\
\[
\chi(X)=10\times (-25)+10\times 5 = -200.
\]
This is exactly the Euler number of a smooth quintic \cy threefold.\\\\
It is instructive to understand the singular point set of the fibration map $f$, namely the points in $X_{\psi}$ where tangent maps of $f$ are not surjective. The singular point set of a type $I_5$ fiber is five $S^1$'s. The singular point set of a type $II_{5\times 5}$ fiber is the graph in figure 2. The singular point set of a type $III_5$ fiber is five points. Together they form a Riemann surface. The Riemann surface has $10$ irreducable components $\{\Sigma_{ijk}\}$. $\Sigma_{ijk}$ is fibered over $\Gamma_{ij}^k\cup \Gamma_{ki}^j\cup \Gamma_{jk}^i$ and centered over $P_{ijk}$. Singular fibers are three of $5$-point sets and $\Gamma_{ijk}$. Total Euler number is \\
\[
\chi(\Sigma_{ijk})=3\times 5 - 25 = -10 = 2-2g(\Sigma_{ijk}).
\]
So $g(\Sigma_{ijk})= 6$.\\\\
Actually there is a very explicit description of the local structure of the map $f$ arround a singular point in a type $III_5$ fiber. It was given in the celebrated paper of R. Harvey and H. B. Lawson (\cite{HL}). The construction is as follows.\\
\begin{th}
Let $M_c=f^{-1}(c)$, where $f:\ {\bf C^3}\longrightarrow {\bf R^3}$\\
\[
f_i(z)= |z_1|^2 - |z_i|^2,\ \ i=2, 3;\ \ \ f_1= {\rm Im}(z_1 z_2 z_3).
\]
Then $M_c$ (with the correct orientation) is a special Lagrangian submanifold of ${\bf C^3}$.
\end{th}
\begin{flushright} $\Box$ \end{flushright}
It is pretty obvious why this characterize the singular points in type $III_5$ singular fiber. Since each $M_c$ is invariant under the action of the group\\
\[
T^2=\left(
\begin{array}{ccc}
e^{i\theta_1} &  0 & 0\\
0 & e^{i\theta_2}  & 0\\
0 & 0 & e^{i\theta_3}
\end{array}
\right)\ \ \theta_1+\theta_2+\theta_3=0.
\]
An easy computation will show that the singular point set around origin is exactly the union of the three coordinate axes. This is where three of the irreducible components of the singular point set meet. The singular locus in ${\bf R^3}$ is\\
\[
\{c\in{\bf R^3}| c_2\leq 0, c_3=c_1=0\}\cup \{c\in{\bf R^3}| c_3\leq 0, c_2=c_1=0\}
\]
\[
\cup \{c\in{\bf R^3}| c_2+c_3\geq 0, c_2-c_3=c_1=0\}.
\]
Fiber over origin is of type $III_5$ singularity and fibers over other points in the singular locus are Type $I_5$. An important observation is that $\Sigma_{ijk}$ is exactly equivalent to the genus 6 curve\\
\[
\{[z]\in {\bf CP^4}| z_i^5 + z_j^5 + z_k^5 =0, z_l=0\ \ {\rm for}\ l\in\{1,2,3,4,5\}\backslash\{i,j,k\}\}
\]
Associate to the special lagrangian torus fiberation $f$ there is the Leray spectral sequence, which abuts to $H^\cdot(X_{\psi}, {\bf Q})$ and in which\\
\[
E_2^{p,q} = {\bf H}^p(\partial \Delta; R^qf_*{\bf Q})
\]
This spectral sequence degenerates at $E_2$ term. We can use it to compute cohomology of $X_{\psi}$.\\\\
This spectral sequence was discussed and computed in \cite{GW} for expected SYZ \sl fiberation of generic \cy manifolds. Our situation is far from generic. Fermat type quintics have a lot of symmetries and actually represent singular (orbifold) points in the moduli space of \cy. It is interesting to compute the spectral sequence of the fiberation and compare with the generic situation.\\\\
Since $R^0f_*{\bf Q} ={\bf Q}$, $\partial \Delta \cong S^3$. It is easy to see that the bottom row of the spectral sequence is\\
\[
E_2^{3,q}\cong H^q(S^3,{\bf Q}):\ \ \ {\bf Q},\ 0,\ 0,\ {\bf Q}.
\]
In general, we have the following exact sequence\\
\[
0\longrightarrow K_q \longrightarrow R^qf_*{\bf Q} \longrightarrow i_* E_q \longrightarrow 0
\]
with the second non-trivial map being the attachment map with $E_q=i^* R^qf_*{\bf Q}$ for $i:\ \partial\Delta\backslash \Gamma \longrightarrow \partial\Delta$. Space $B=\partial\Delta$ has a natural filtration\\
\[
B=B_3\supset B_2\supset B_1\supset B_0\supset B_0 = \phi
\]
where $B_2=B_1=\Gamma$, $B_0=\{P_{ij}|i,j\in \{1,2,3,4,5\}\}\cup \{P_{ijk}|i,j,k\in \{1,2,3,4,5\}\}$. Each stratum $S_i=B_i\backslash B_{i-1}$ is a $i$-manifold. We denote $V_i=B\backslash B_i$. $E_q$ is a local system on $V_1=\partial\Delta\backslash \Gamma$ and $K_q$ is supported on $B_1=\Gamma$. Each $R^qf_*{\bf Q}$ is a locally constant sheaf. Its germ over $P\in B$ is\\
\[
(R^qf_*{\bf Q})_P\cong H^q(f^{-1}(P), {\bf Q}) \cong (H_q(f^{-1}(P), {\bf Q}))^\vee
\]
$E_3$ is obviously a constant sheaf, so $i_*E_3\cong {\bf Q}$. According to our geometric description of the singular fibers, $R^3f_*{\bf Q}$ is of dimension $1$ over $S_3$, dimension $5$ over $S_1$, dimension $5$ over $P_{ij}$ and dimension $25$ over $P_{ijk}$. Therefore $K_3$ is supported in $B_1=\Gamma$ and is of dimension $4$ over $S_1$, dimension $4$ over $P_{ij}$ and dimension $24$ over $P_{ijk}$. Clearly only possible non-trivial cohomologies for $K_3$ are $H^0(K_3,{\bf Q})$ and $H^1(K_3,{\bf Q})$.\\\\
\begin{prop}
\[
H^1(K_3,{\bf Q})=0
\]
\end{prop}
{\bf Proof:} We will compute the {\v C}ech cohomology. There is a natural map\\
\[
\pi_{ijk}:\ (K_3)_{P_{ijk}} \longrightarrow  (K_3)_{\Gamma^k_{ij}}\oplus (K_3)_{\Gamma^i_{jk}}\oplus (K_3)_{\Gamma^j_{ki}}
\]
We will show that $\pi_{ijk}$ is surjective, which easily implies that $H^1(K_3,{\bf Q})=0$. Observe that $\pi_{ijk}$ is induced from\\
\[
{\tilde \pi}_{ijk}: (R^3f_*{\bf Q})_{P_{ijk}} \longrightarrow  (R^3f_*{\bf Q})_{\Gamma^k_{ij}}\oplus (R^3f_*{\bf Q})_{\Gamma^i_{jk}}\oplus (K_3)_{\Gamma^j_{ki}}
\]
Choose generators $\{x_{ij}|i,j\in \{1,2,3,4,5\}\}$ of $(R^3f_*{\bf Q})_{P_{ijk}}$, $\{u_i\}_{i=1}^5$ of $(R^3f_*{\bf Q})_{\Gamma^k_{ij}}$, $\{v_i\}_{i=1}^5$ of $(R^3f_*{\bf Q})_{\Gamma^j_{ki}}$ and $\{w_i\}_{i=1}^5$ of $(R^3f_*{\bf Q})_{\Gamma^k_{ij}}$ (representing geometric cycles). If we choose correctly, we should have\\
\[
{\tilde \pi}_{ijk}(x_{lm})=(u_l, v_m, w_n),\ \ n \equiv -l-m\ (mod\ 5).
\]
By simple linear algebra, we can see that\\
\[
\left( \sum^5_{i=1}a_l u_l,\ \sum^5_{i=1}b_l v_l,\ \sum^5_{i=1}c_l w_l\right) \in {\rm Im}({\tilde \pi}_{ijk})
\]
if and only if
\[
\sum^5_{i=1}a_l =\sum^5_{i=1}b_l =\sum^5_{i=1}c_l;
\]
and
\[
\left( \sum^5_{i=1}a_l u_l,\ \sum^5_{i=1}b_l v_l,\ \sum^5_{i=1}c_l w_l\right) \in {\rm Im}(\pi_{ijk})
\]
if and only if
\[
\sum^5_{i=1}a_l =\sum^5_{i=1}b_l =\sum^5_{i=1}c_l =0.
\]
Therefore $\pi_{ijk}$ is surjective.
\begin{flushright} $\Box$ \end{flushright}
\begin{prop}
\[
h^0(K_3)=160.
\]
\end{prop}
{\bf Proof:}  We will compute Euler number $e(K_3)=h^0(K_3)-h^1(K_3)$. We compute in the chain level. Let $c^i=\dim(C^i)$ be the dimension of the space of dimension $i$ cochain. It is easy to see\\
\[
c^0=10\times 24 + 10 \times4 =280,\ \ c^1=30\times4=120,\ \ e(K_3)=c^0-c^1=160.
\]
Then $h^1(K_3)=0$ from last proposition gives us that $h^0(K_3)=160.$
\begin{flushright} $\Box$ \end{flushright}
Now use the long exact sequence associated to\\
\[
0\longrightarrow K_3 \longrightarrow R^3f_*{\bf Q} \longrightarrow i_* E_3 \longrightarrow 0
\]
we get
\begin{prop}
\[
h^0(R^3f_*{\bf Q})=161,\ h^1(R^3f_*{\bf Q})=h^2(R^3f_*{\bf Q})=0,\ h^3(R^3f_*{\bf Q})=1.
\]
\end{prop}
\begin{flushright} $\Box$ \end{flushright}
To compute the cohomology for $R^1f_*{\bf Q}$, notice that $R^1f_*{\bf Q}=i_*E_1$ for local system $E_1$. It is easy to see that the monodromy of the local system $E_1$ generate $SL_3({\bf Z})$. So $h^0(R^1f_*{\bf Q})=0$. From this we will have\\
\begin{prop}
\[
h^3(R^1f_*{\bf Q})=h^3(R^2f_*{\bf Q})=0.
\]
\end{prop}
{\bf Proof:} By Poincare duality for general fibers, $E_2=E_1^\vee$. So\\
\[
R^2f_*{\bf Q}\cong {\cal D}i_*E_1 \ \ \ {\rm in}\ V_2
\]
where ${\cal D}E$ is the dual of $E$ in the sense of Verdier. Then by Verdier duality\\
\[
H^3(R^2f_*{\bf Q})\cong H^3({\cal D}i_*E_1)\cong H^0(i_*E_1)^\vee=0.
\]
The space $B$ have a natural symmetry $s$ which respects the filtration. It is a piecewise linear map satisfying\\
\begin{equation}
\label{ca}
s^{-1}=s:\ B \longrightarrow B;\ s(P_i)=P_{\overline{i}},\ s(P_{ij})=P_{\overline{ij}}
\end{equation}
where $\overline{S}$ indicate compliment of $S\in \{1,2,3,4,5\}$. It is easy to check by looking at the monodromy that $s^*E_1\cong E_2$. (Later we will see that this fact is related to the mirror symmetry construction of Fermat type mirror construction.) Therefore\\
\[
i_*E_1\cong s^*R^2f_*{\bf Q}\ \ \ {\rm in}\ V_2
\]
and
\[
H^3(R^1f_*{\bf Q})=H^3(i_*E_1)\cong H^3(R^2f_*{\bf Q})=0.
\]
\begin{flushright} $\Box$ \end{flushright}
To compute $H^1(R^1f_*{\bf Q})$, we need to use intersection cohomology. Since $R^1f_*{\bf Q}=i_*E_1$, We have $H^1(R^1f_*{\bf Q})=I_0 H^1(B,E_1)$. Where $I_p H^i(B,E)$ denote the $i$-th intersection cohomology of $B$ with coefficient in the local system $E$ and perversity $p$.\\
\begin{prop}
\label{bd}
\[
h^1(R^1f_*{\bf Q})=1.
\]
\end{prop}
{\bf Proof:} By Poincare duality and $E_2=E_1^\vee$, we have\\
\[
H^1(R^1f_*{\bf Q})=I_0 H^1(B,E_1)=I_t H^2(B,E_2)^\vee
\]
where $t$ denote the top perversity. Intersection homology is related to intersection cohomology as\\
\[
I_t H^2(B,E_2)=I_t H_1(B,E_2).
\]
So we only need to compute $I_t H_1(B,E_2)$, where for $P\in V_1$ we have $(E_2)_P\cong H_1(f^{-1}(P),{\bf Q})$. Consider the first barycentric subdivision of the standard triangulation of $B$. The 1-simpleces allowed to compute $I_t H_1(B,E_2)$ are the 1-simpleces $\Gamma^i_{jkl}$ connecting $P_i$ and $P_{ijkl}$.\\
\begin{figure}[h]
\begin{center}
\leavevmode
\hbox{%
\epsfxsize=5.5in
\epsffile{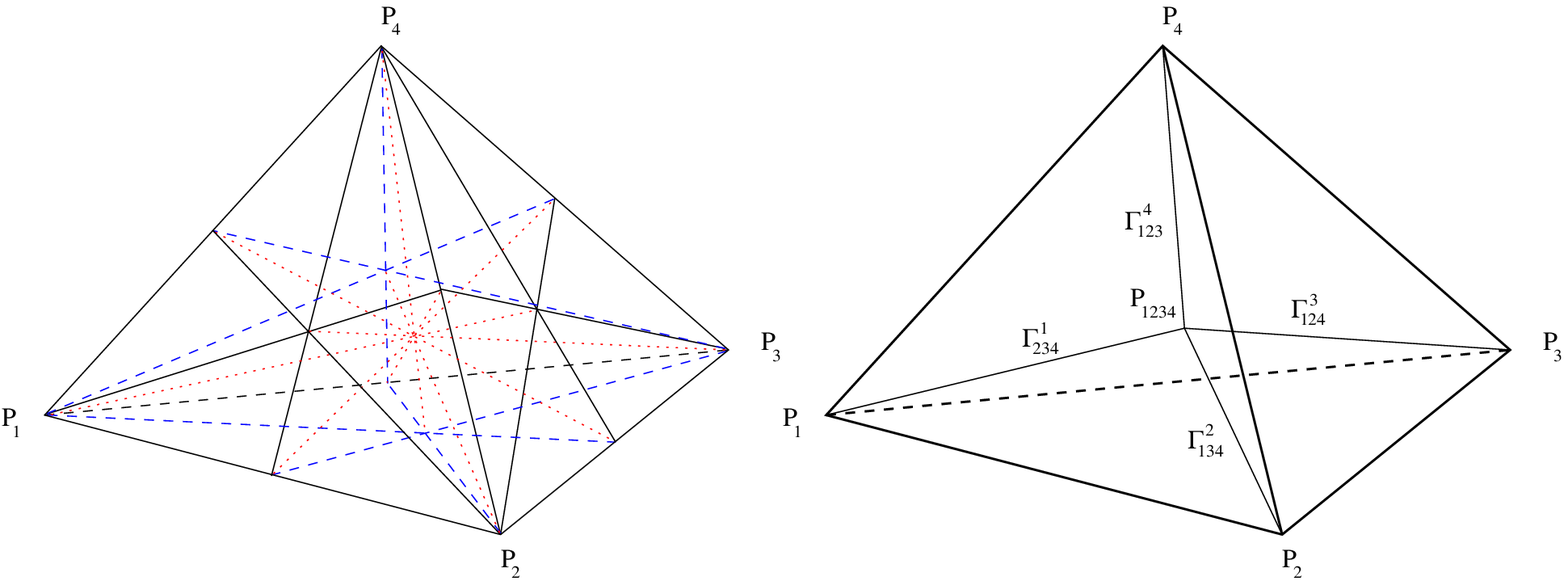}}
\end{center}
\caption{}
\end{figure}\\
We can introduce the following 1-chain\\
\[
L=\gamma_j^i \cdot \Gamma^i_{\overline{ij}}
\]
By (\ref{bc}), it is easy to see that $\partial L=0$. Actually this $L$ generate $I_t H_1(B,E_2)$, namely\\
\[
h^1(R^1f_*{\bf Q}) = \dim (I_t H_1(B,E_2))=1
\]
\begin{flushright} $\Box$ \end{flushright}
\begin{prop}
\label{be}
\[
h^2(R^1f_*{\bf Q})=1.
\]
\end{prop}
{\bf Proof:} We will show this by proving that $\chi (R^1f_*{\bf Q})=0$, which clearly implies this proposition by results in proposition \ref{bd} and \ref{be}. The Euler number can be computed in two ways, either by straightforward {\v C}ech cohomology computation or by the identification\\
\[
H^i(R^1f_*{\bf Q}) \cong I_0 H^i(B,E_1)\cong I_0 H_{3-i}(B,E_1).
\]
We will use the latter approach, which is more elegant. To get the number right, the key point is to choose the right triangulation, which respects the filtration of $B$. The previous trangulation as indicated in figure 4 is not good enough. The problem is that the 0-simplices $P_{ij},P_{ijk}\in B_0$ will span a 1-simplex $\Gamma^k_{ij}$ which is in $B_1$ but not in $B_0$. The right triangulation which respects the filtration of $B$ can be achieved by taking the baricentric subdivision of $\Gamma$ and then naturally extend to a subdivision of the previous triangulation of $B$. Practically, each old 3-simplex is divided into two new 3-simleces by a new 2-simplex as indicated in the following picture.\\
\begin{figure}[h]
\begin{center}
\leavevmode
\hbox{%
\epsfxsize=5.5in
\epsffile{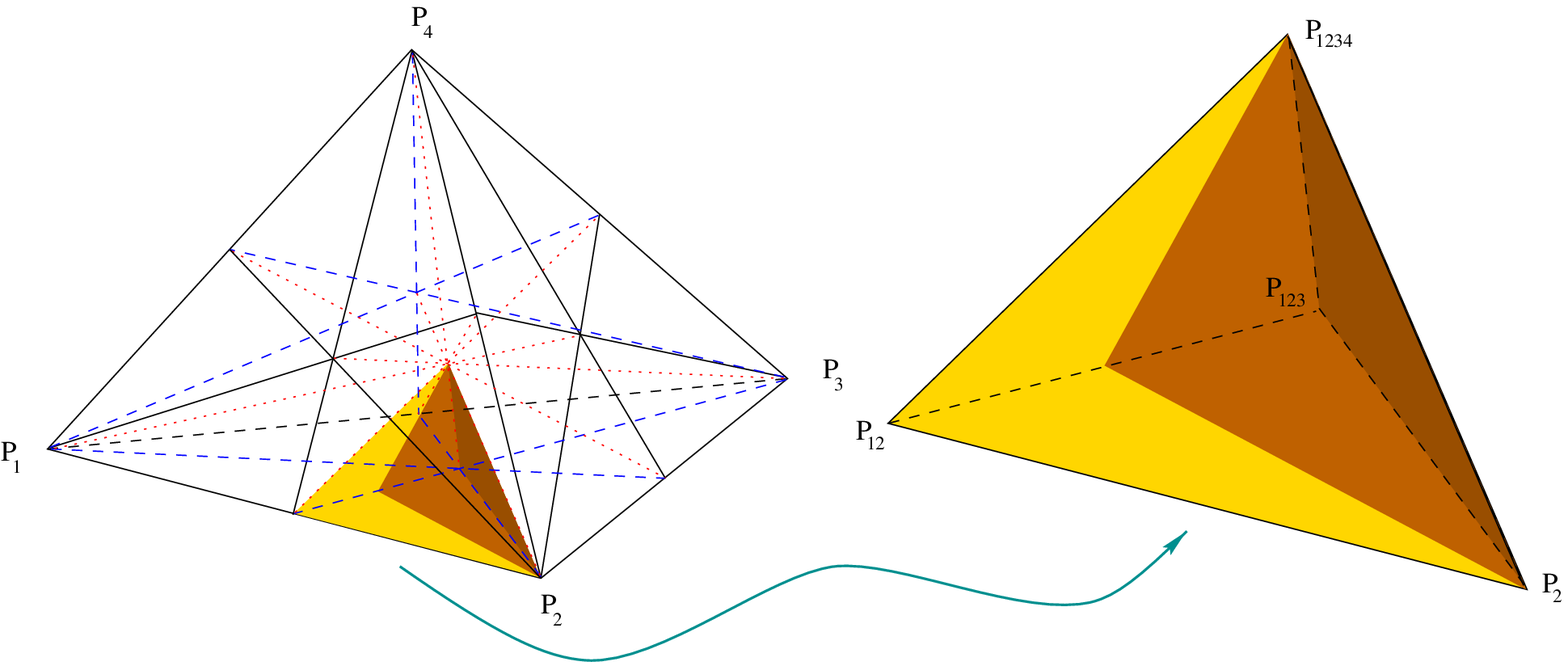}}
\end{center}
\caption{}
\end{figure}\\
It is easy to see that\\
\[
I^0C_0=\bigoplus_{i=1}^5(E_1)_{P_{i}}\cdot P_{i} \oplus (E_1)_{P_{\overline{i}}}\cdot P_{\overline{i}}
\]
\[
I^0C_1=\bigoplus_{\tiny{\begin{array}{c} i,j=1\\ (i\not= j)\end{array}}}^5(E_1)_{<P_{i},P_{\overline{j}}>}\cdot <P_{i},P_{\overline{j}}>,
\]
and $c_0=10\times 3=30$, $c_1=20\times 3=60$.\\\\
For the $I^0C_2$, it is easy to see that only the newly added 2-simplices $<P_{i}, P_{\{ij\}\{ijk\}}, P_{ijkl}>$ will be relavent. Since the boundary of the intersection 2-chain should not intersect $\Gamma$. An intersection 2-chain containing $<P_{i}, P_{\{ij\}\{ijk\}}, P_{ijkl}>$ should also contain all other newly added 2-simplices that contain $P_{\{ij\}\{ijk\}}$, namely should contain\\
\[
\Delta^k_{ij}=\bigcup\{{\rm 2-simplex\  with}\  P_{\{ij\}\{ijk\}}\ {\rm as\ one\ vertex}\}.
\]
Therefore\\
\[
I^0C_2=\bigoplus_{\Gamma_{ij}^k}(E_1)^{T^k_{ij}}_{\Delta^k_{ij}}\cdot \Delta^k_{ij}
\]
where the upper index ${T^k_{ij}}$ indicates the ${T^k_{ij}}$ invariant piece of $E_1$. So $c_2=30\times2=60$.\\\\
Since an intersection 2-chain is not supposed to contain $P_{ij}$ or $P_{ijk}$. The boundary of an intersection 3-chain also should not contain $P_{ij}$ or $P_{ijk}$. Hence an intersection 3-chain containing $<P_{i}, P_{ij}, P_{\{ij\}\{ijk\}}, P_{ijkl}>$ has to contain\\
\[
\Delta_{ij}=\bigcup\{{\rm 3-simplex\  with}\  P_{ij}\ {\rm as\ one\ vertex}\}.
\]
An intersection 3-chain containing $<P_{i}, P_{ijk}, P_{\{ij\}\{ijk\}}, P_{ijkl}>$ has to contain\\
\[
\Delta_{ij}=\bigcup\{{\rm 3-simplex\  with}\  P_{ijk}\ {\rm as\ one\ vertex}\}.
\]
And\\
\[
I^0C_3=\left(\bigoplus_{ij} (E_1)^{G_{ij}}_{\Delta_{ij}}\cdot \Delta_{ij}\right) \bigoplus \left(\bigoplus_{ijk} (E_1)^{G_{ijk}}_{\Delta_{ijk}}\cdot \Delta_{ijk}\right),
\]
where $G_{ij}$ is the monodromy group around $P_{ij}$ and $G_{ijk}$ is the monodromy group around $P_{ijk}$. And $c_3=10\times 1+10\times2=30$. Therefore\\
\[
\chi(R^1f_*{\bf Q})=-\chi(I^0C.)=-\sum_{i=0}^3 (-1)^i c_i = -(30-60+60-30)=0.
\]
\begin{flushright} $\Box$ \end{flushright}
Now the only thing we do not yet know is $h^1(R^2f_*{\bf Q})$.\\
\begin{prop}
\[
h^1(R^2f_*{\bf Q})=41.
\]
\end{prop}
{\bf Proof:} Clearly we only need to compute $h^1(K_2)$. Since the germ of $K_2$ at $P_{ij}$ only have zero section. It is clear that $h^0(K_2)=0$. It is now reduced to compute $\chi(K_2)$. We will simply count the {\v C}ech cochains.\\
\[
c^1=(5-1)\times 30=120,\ \ c^0=(10-2)\times 10 = 80.
\]
\[
\chi(K_2)=c^0-c^1 = 80-120=-40.
\]
So $h^1(K_2)=40$ and $h^1(R^2f_*{\bf Q})=41$.
\begin{flushright} $\Box$ \end{flushright}
Now we have the whole $E_2$ term of the spectral sequence.\\
\begin{th}
$E_2$ term of the Leray spectral sequence of fibration $F$ is\\
\[
\begin{array}{cccc}
161 & 0 & 0 & 1\\
0 & 41 & 1 & 0\\
0 & 1 & 1 &0\\
1 & 0 & 0 & 1
\end{array}
\]
\end{th}
\begin{flushright} $\Box$ \end{flushright}
One can see that the $E_2$ term of our situation is quite different from the $E_2$ term of generic \cy situation as indicated in \cite{GW}.\\\\
Note: Our discussion so far is of conjectural nature. Since we can not construct \sl fibration, we can only guess its possible structure based on our monodromy discussion.\\

\se{Comparison}
The \l fibration we constructed in Section 3 is different from the expected \sl fibration structure described in the previous two sections. Yet they are very closely related. The singular locus $\Gamma$ of the expected \sl fibration is of codimension 2 (an one dimensional graph) and each singular fiber has singularity of codimension 2. The singular locus $\tilde{\Gamma}$ of the \l fibration we constructed is of codimension 1 and most singular fibers have singularity of codimension 3. On the other hand, $\tilde{\Gamma}$ is just a fattened version of $\Gamma$. $S^3\backslash\tilde{\Gamma}$ and $S^3\backslash\Gamma$ have the identical fundamental group. We can compare monodromies of the two fibrations.\\\\
Since our \l fibration on $X_{\psi}$ is a deformation of the standard \sl fibration on $X_{\infty}$, it satisfies the principle assumption (1) that we based to compute monodromy. For this reason, it is not hard to see that above computation of monodromy naturally apply to our \l fibration. Therefore\\
\begin{th}
Our \l fibration have the same monodromy as the monodromy computed for the expected \sl fibration.
\end{th}
\begin{flushright} $\Box$ \end{flushright}
The singularities of the two fibrations are also closely related. For one thing, the singular point set of our construction in $X_{\psi}$ is $\Sigma$ (the union of 10 genus 6 curves), and the singular point set of the expected \sl fiberation is naturally equivalent to $\Sigma$. One way to see this is to notice that we can easily construct a homotopy contraction $p: \tilde{\Gamma} \rightarrow \Gamma$ that is the homotopy inverse of the natural injection $i: \Gamma \rightarrow \tilde{\Gamma}$ that satisfy $p\cdot i= id$ (for instance map point in $\tilde{\Gamma}$ to the nearest point in $\Gamma$). Let $q$ denote the composition of projection of $\Sigma$ to $\tilde{\Gamma}$ and $p$. Then $q$ map $\Sigma$ to $\Gamma$. It is easy to observe that the inverse image of a point $P$ in $\Gamma$ with respect to $q$ is exactly corresponding to singular set of the expected singular \sl fiber over $P$ as described in the previous section. With this fact in mind, we intend to modify our \l torus fibration into the expected topological shape.\\\\
If one is willing to sacrifies \l property, it is not hard to deform the fibration we constructed into a non-\l smooth fiberation with the same topological shape as we expected of \sl fibration. One way to do this is to notice that there is a natural horizontal foliation in ${\bf P^2}$ where one of the component of $\Sigma$ is located. Then conceivably, one can smoothly deform part of 2-torus above $\tilde{\Gamma}$ that intersect $\Sigma$ along horizontal direction away from $\Gamma$ in direction indicated by map $p$. So that eventually only those 2-torus above $\Gamma$ will intersect $\Sigma$ at one dimensional graphs that is the inverse image under $q$ of the corresponding point. Then we run the flow of $V$, we will get a non-\l torus fibration with topological structure that is the same as expected \sl fibration. This will give an alternative simple proof of the result in \cite{Z} that gave a smooth torus non-\l fibration. Since our major concern is to construct smooth \l torus fiberation, we will not get into great details in this direction.\\\\
However to get a \l torus fibration of the right shape is much trickier. We will address this problem in the sequel to this  paper(\cite{lag2}).\\\\
Another important observation is that for an $n$-torus \l fibration such that the fibration map is $C^{\infty}$ (actually $C^2$ is enough), there is a natural action of ${\bf R^n}$ on each fiber (smooth or singular). Look at our construction, one can see that there are no action of ${\bf R^3}$ on most of our singular fibers, which seems to be a contradiction. It turns out that fibration map of our construction is not $C^{\infty}$---it is typically only piecewise smooth and globally merely $C^{0,1}$ (Lipschitz). This reveals a major difference between a non-\l fibration and a \l fibration. If one does not care keeping the \l property, the map usually can be smoothed by a small perturbation. Above ovservation shows that in general a \l fiberation can not be deformed to a smooth \l fiberation by small perturbation. In (\cite{lag2}) we will explore this issue further and try to get a smooth \l fibration. We will also discuss general construction of \l torus fibration for general quintic \cy hypersurfaces in \cite{lag3} and \cy hypersurfaces in more general toric varieties in \cite{tor}.\\\\

\se{The Mirror construction and dual polyhedra}
It is well known that the mirror familly $\{Y_\psi\}$ of the family $\{X_\psi\}$ can be given by $Y_\psi=X_\psi/G$ with the induced complex structure and the \cy orbifold metric. $G=<g_0, g_1, g_2, g_3>$ is the symmetry group of $X_\psi$. $G$ acts on ${\bf CP^4}$ as follows\\
\begin{eqnarray*}
g_0&=&(\xi, 1, 1, 1, \xi^{-1}),\\
g_1&=&(1, \xi, 1, 1, \xi^{-1}),\\
g_2&=&(1, 1, \xi, 1, \xi^{-1}),\\
g_3&=&(1, 1, 1, \xi, \xi^{-1}),
\end{eqnarray*}
where $\xi$ is a fifth root of unity. By analyzing the action of $G$ on $X_\psi$, it should be clear that $G\cong {\bf Z}_5^3$ will map every special Lagrangian fiber to itself, and is conjugate to the standard action of ${\bf Z}_5^3$ on $T^3$. The special Lagrangian fibration $f:\ X_\psi\rightarrow B$ will naturally induce special Lagrangian fibration $\hat{f}:\ Y_\psi\rightarrow B$. For any $p\in B$, $\hat{f}^{-1}(p)=f^{-1}(p)/G$, especially $\hat{f}^{-1}$ can be identified with $f^{-1}(p)$ quite canonically.\\\\
SYZ construction predict that the $T^3$ fibrations $f$ and $\hat{f}$ when restrict to $V_1\subset B$ should be dual to each other. But the above discussion implies that the two fibrations can be natually identified to each other instead of dual, especially the monodromy for $\hat{f}$ is the same as the monodromy for $f$. The key point is that the mirror map induce nontrivial identification $s:\ B\rightarrow B$ of base $B$ as defined in (\ref{ca}). Fibrations $\hat{f}$ and $B\circ f$ are actually dual to each other when restricted to $V_1\subset B$. This corresponds to the fact that the monodromy around $P_{ij}$ is dual to the monodromy around $P_{\overline{ij}}$. The reason why this is the correct interpretation is that the two fibrations should actually be written as\\
\[
f:\ X_\psi\longrightarrow \Delta,\ \ \ \hat{f}:\ Y_\psi\longrightarrow \hat{\Delta}
\]
where $\hat{\Delta}$ is the dual polyhedron of $\Delta$. The only canonical map from $\Delta$ to $\hat{\Delta}$ is the one corresponding to map $s$.\\\\
Now let us analyze the singular fibers of $\hat{f}$. For the type $I_5$ singular fiber $I_5\times S^1$, one ${\bf Z}_5$ will rotate $S^1$, another ${\bf Z}_5$ will rotate each $S^2$ in $I_5$ while keeping the nodal points fixed, the ${\bf Z}_5$ left will permute the five $S^2$'s in $I_5$. The quotient will natually be $I_1\times S^1$, we call it type $I$. A type $I$ fiber has Euler number $0$.\\\\
For the type $II_{5\times 5}$ singular fiber over $P_{ijk}$, two ${\bf Z}_5$'s will act on $T^2\cong T^3/{\gamma^m_l}$ as indicated by the two arrows in the following picture, which move a sixgon along the arrows. The ${\bf Z}_5$ left will act on $\gamma^m_l$ in standard way, which do not affect topology. Upon the action, edges of a fundamental reigon  in the $T^2$ as indicated as darker area in the picture will be identified as indicated in the picture and turn into a $T^2$, with image of the graph $\Gamma_{ijk}$ turn into $\hat{\Gamma}_{ijk}$, topologically, a contractor of ``a pair of pants''. The resulting singular fiber can be identified as $S^1\times T^2$ with $S^1$'s over $\hat{\Gamma}_{ijk}$ collapsed to points. We will call it type $II$ singular fiber. A type $II$ fiber has Euler number $1$.\\
\begin{figure}[h]
\begin{center}
\leavevmode
\hbox{%
\epsfxsize=5in
\epsffile{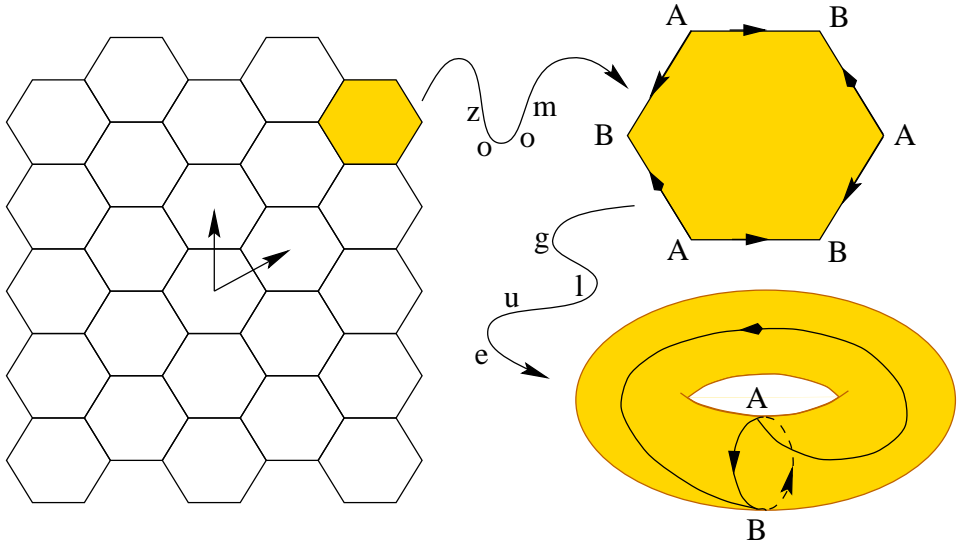}}
\end{center}
\caption{Type $II$ fiber}
\end{figure}\\
The quotients of the type $III_5$ singular fibers are very easy to describe. Two ${\bf Z}_5$'s will act on $T^2$ in standard way, topologically cause no change. The ${\bf Z}_5$ left will rotate among the chain of the five suspensions of $T^2$. The resulting singular fiber will simply be a suspension of $T^2$ with the two pole points identified as indicated in the following picture. We will call it type $III$ fiber. A type $III$ fiber has Euler number $-1$.\\
\begin{figure}[h]
\begin{center}
\leavevmode
\hbox{%
\epsfxsize=5in
\epsffile{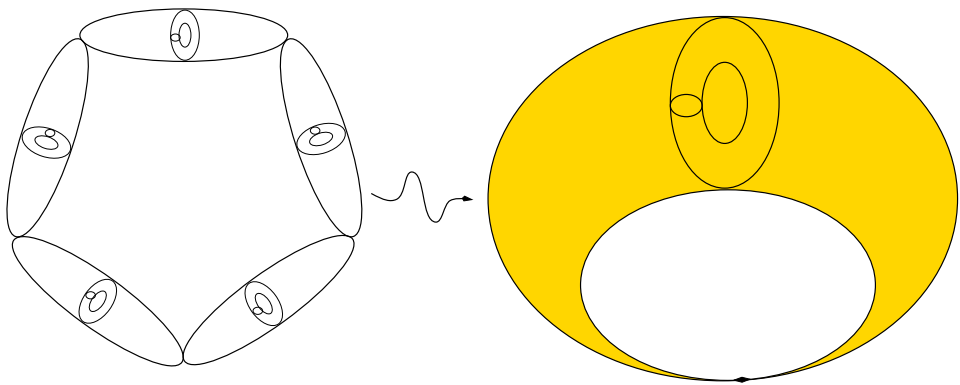}}
\end{center}
\caption{Type $III$ fiber}
\end{figure}\\
It is interesting to see what is the image of the singular point sets\\
\[
\Sigma_{ijk}=\{[z]\in {\bf CP^4}| z_i^5 + z_j^5 + z_k^5 =0, z_l=0\ \ {\rm for}\ l\in\{1,2,3,4,5\}\backslash\{i,j,k\}\}
\]
under the quotient map. It is easy to see that $G$ map $\Sigma_{ijk}$ to itself. Two of the ${\bf Z}_5$'s act as identity, and the ${\bf Z}_5$ left acts freely on $\Sigma_{ijk}$. So $\Sigma_{ijk}$ is a 5-sheet cover of the quotient $\hat{\Sigma}_{ijk}$. Since\\
\[
5\cdot \chi(\hat{\Sigma}_{ijk})=\chi(\Sigma_{ijk})=-10
\]
we have\\
\[
\chi(\hat{\Sigma}_{ijk})=-2,\ \ g(\hat{\Sigma}_{ijk})=0.
\]
Namely $\hat{\Sigma}_{ijk}$ is a Riemann sphere. Their union is exactly the singular point set of the quotient $Y_{\psi}$.\\\\
Now we can try to compute the Leray sqectral sequence of the fibration $\hat{f}$. By analyzing the singular fibers described at above, we can esaily see that\\
\[
R^0 \hat{f}_*{\bf Q}\cong R^3 \hat{f}_*{\bf Q}\cong {\bf Q}
\]
and
\[
R^1 \hat{f}_*{\bf Q}\cong (R^2 \hat{f}_*{\bf Q})^\vee \cong i_* E_1
\]
These combine with computation from last section, give us\\
\begin{prop}
$E_2$ term of the Leray spectral sequence of fibration $\hat{f}$ is\\
\[
\begin{array}{cccc}
1 & 0 & 0 & 1\\
0 & 1 & 1 & 0\\
0 & 1 & 1 & 0\\
1 & 0 & 0 & 1
\end{array}
\]
\end{prop}
\begin{flushright} $\Box$ \end{flushright}
What follows is a discussion of the well known crepant resolution of $Y_{\psi}$. For completeness of mirror picture, we include it here. We need some degression on toric varieties. Let $X=Spec({\bf C}[C])$ be an affine toric variety, where
\[
C=<v^1, v^2, \cdots, v^n>_{\bf Q^+}\cap M \subset M=<e^1, e^2, \cdots, e^n>_{\bf Z}
\]
is a cone in lattice $M$. Following convention, let $N$ be the dual lattice of $M$ and ${\check C}$ be the dual cone of $C$. Suppose
\[
{\check C}=<v_1, v_2, \cdots, v_n>_{\bf Q^+}\cap N \subset N=<e_1, e_2, \cdots, e_n>_{\bf Z}.
\]
It is well known that a birational modification $X'$ of $X$ corresponds to subdividing ${\check C}$ into union of smaller cones (we require each new cone to be generated by $n$ elements), one dimensional edges of these cones correspond to toric divisors in $X'$. Let $D_v$ denote the toric divisor corresponding to the 1-dimensional cone generated by $v$. We will use $X_0$ to denote the smooth part of $X$. It is obvious that the singular part of $X$ is of codimension greater or equal to two. We need the following result.\\
\begin{prop}
There is a canonical nowhere vanishing holomorphic n-form $\Omega$ unique up to constant on $X_0$. Assume $X'$ is a birational modification of $X$, then $\Omega$ can be extend to $X_0$ as a nowhere vanishing holomorphic n-form if and only if the generater of any of the 1-dimension subcone coming from the subdivision is a convex combination of $ v_1, v_2, \cdots, v_n$.
\end{prop}
{\bf Proof:} It is well know that there is a canonical holomorphic n-form $\Omega_0$ on the $n$-torus defined as
\[
\Omega_0=\frac{{\rm d}e^1}{e^1}\wedge\frac{{\rm d}e^2}{e^2}\wedge\cdots \wedge\frac{{\rm d}e^n}{e^n}
\]
For any primary $v\in N$, we can get $D_v=Spec({\bf C}[v^\bot])$. $\Omega_0$ will have pole of order 1 around $D_v$. To cancel the pole, we need to multiply a function, which vanish to order one at $D_v$. These kind of function will be in\\
\[
v^\vee =\{w\in M|<v, w>=1\}
\]
Suppose that\\
\[
\bigcap _{i=1}^n v_i^\vee =\{w\}
\]
then $\Omega=w\cdot\Omega_0$ is exactly what we need.\\\\
Suppose that $v$ is the generater of an 1-dimension subcone coming from the subdivision corresponding to a birational modification $X'$. Then in order to extend $\Omega=w\cdot\Omega_0$ to $D_v\subset X'$, $v$ must satisfy $<v, w>=1$, which is the same as that $v$ is a convex combination of $ v_1, v_2, \cdots, v_n$.
\begin{flushright} $\Box$ \end{flushright}
{\bf Remark:} The above proposition must be quite well known in toric geometry. But we do not know exactly where it was located.\\\\
$Y_{\psi}$ is a singular \cy manifold. It has singularity at $\tilde{P}_{ij}$ and $\hat{\Sigma}_{ijk}$, where $\tilde{P}_{ij}$ is the singular point of the fiber over $P_{ij}$. Around $\tilde{P}_{ij}$, $Y_{\psi}$ has quotient singularity. It can be modeled by affine toric variety $X=Spec({\bf C}[C])$, where\\
\[
C=<5e^1, 5e^2, 5e^3>_{\bf Q^+}\cap M \subset M=<5e^1, 5e^2, 5e^3, \sum_{i=1}^3 e^i>_{\bf Z}
\]
Here $\{e^i\}$ denote the standard base of ${\bf Z}^3$. We would like to get a crepent resolution of the toric variety $X=Spec({\bf C}[C])$, for which we have to subdivide the dual cone\\
\[
{\check C}=<e_1, e_2, e_3>_{\bf Q^+}\cap N \subset N=M^\vee,
\]
here $\{e_i\}$ denote the standard base of the dual ${\bf Z}^3$. By above discussion, we only need to look at the plane passing through the three points $e_1, e_2, e_3$. ${\check C}$ will cut off a triangle with this three points as vertices in the plane. All the lattice points in $N$ that is in this triangle are\\
\[
v_{ijk}=\frac{i}{5}\cdot e_1+\frac{j}{5}\cdot e_2 +\frac{k}{5}\cdot e_3,\ \ i,j,k\geq 0\ \ i+j+k=5
\]
We will divide this triangle as indicated in following picture. This division will result in a subdivision of ${\check C}$ which gives us a desired crepent resolution.\\
\begin{figure}[h]
\begin{center}
\leavevmode
\hbox{%
\epsfxsize=5in
\epsffile{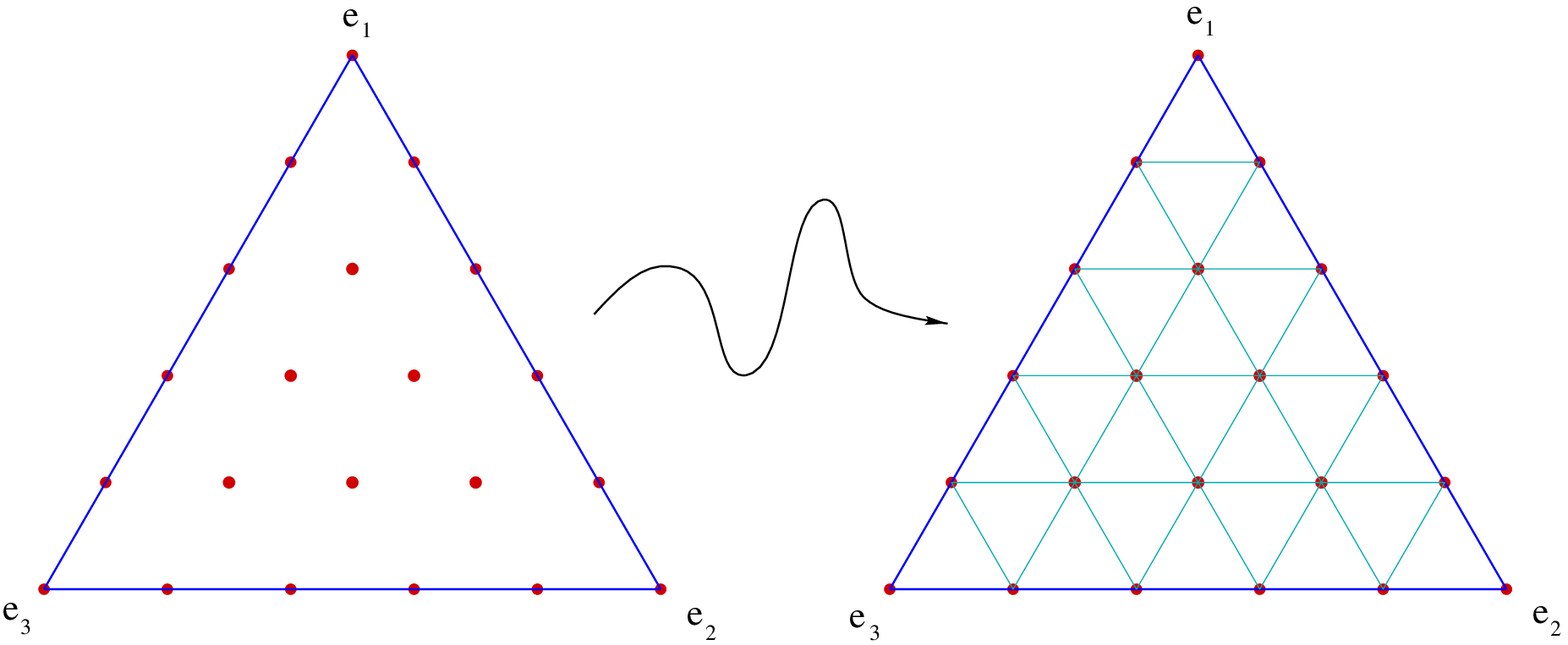}}
\end{center}
\caption{}
\end{figure}\\
There are exactly three curves of singularity $\hat{\Sigma}_{ijk}$ for $k\in \overline{\{ij\}}$ that pass through the singular point $\tilde{P}_{ij}$. In the picture, they correspond to the interiors of the three edges of the first triangle, while the singular point correspond to the interior of the triangle. After the birational modification, we will have newly added exceptional divisors that correspond to dots other than $e_1,e_2,e_3$. The ones in the interiors of the three edges correspond to exceptional divisors coming from the three curves of singularity. The ones in the interior of the first triangle correspond to exceptional divisors coming from the singular point. So each curve of singularity will give rise to 4 exceptional divisors and each singular point will give rise to 6 exceptional divisors. There are 10 curves of singularity and 10 singular points in $Y_\psi$. So all together we will have $4\times 10 + 6\times 10= 100$ of exceptional divisors in the smooth \cy $\tilde{Y}_\psi$ we get as the crepant resolution of $Y_\psi$. This verify the fact\\
\begin{prop}
\[
h^0(\tilde{Y}_\psi)=h^6(\tilde{Y}_\psi)=1, \ h^1(\tilde{Y}_\psi)=h^5(\tilde{Y}_\psi)=0,\ h^2(\tilde{Y}_\psi)=h^4(\tilde{Y}_\psi)=101,\ h^3(\tilde{Y}_\psi)=4.
\]
\end{prop}
\begin{flushright} $\Box$ \end{flushright}

{\bf Acknowledgement:} I would like to thank Qin Jing for many very stimulating discussions during the course of my work, and helpful suggestions while carefully reading my early draft. I would also like to thank Prof. Yau for his constant encouragement.\\\\


\begin{thebibliography}{11}

\bibitem{Bat}
Batyrev, V.~V. ``Dual polyhedra and mirror symmetry for {C}alabi-{Y}au
  hypersurfaces in toric varieties", J. Algebraic Geom. {\bf 3} (1994),
  493--535.
%
\bibitem{Borel}
Borel, A. et al,
{\it Intersection Cohomology},
Birkhauser 1984.
%
\bibitem{Can}
Candelas, P., de la Ossa, X.C., Green, P., Parkes, L.,
"A Pair of \cy Manifolds as an Exactly Soluble Superconformal Theory",
in {\it Essays on Mirror Symmetry}, edited by S.-T. Yau.
%
\bibitem{Gross1}
Gross, M.,
"Special Lagrangian Fibration I: Topology",
alg-geom 9710006.
%
\bibitem{Gross2}
Gross, M.,
"Special Lagrangian Fibration II: Geometry",
alg-geom 9809072.
%
\bibitem{GW}
Gross, M. and Wilson, P.M.H.,
"Mirror Symmetry via 3-torus for a class of \cy Threefolds",
Math. Ann.  309  (1997),  no. 3, 505--531.
%
\bibitem{HL}
Harvey, R. and Lawson, H.B.,
"Calibrated Geometries",
{\it Acta Math.} 148 (1982), 47-157.
%
\bibitem{H}
Hitchin, N.,
"The Moduli Space of Special Lagrangian Submanifolds",
dg-ga 9711002
%
\bibitem{lag2}
Ruan, W.-D.,
``Lagrangian torus fibration of quintic \cy hypersurfaces II: Technical results on gradient flow construction",
Journal of Symplectic Geometry, Volume 1 (2002), Issue 3, 435-521.
%
\bibitem{lag3}
Ruan, W.-D.,
``Lagrangian torus fibration of quintic \cy hypersurfaces III: Symplectic topological SYZ mirror construction for general quintics",
Journal of Differential Geometry, Volume 63 (2003), 171--229.
%
\bibitem{N}
Ruan, W.-D.,
``Newton polygons and string diagrams",
math.DG/0011012.
%
\bibitem{tor}
Ruan, W.-D.,
``Lagrangian torus fibration and mirror symmetry of \cy hypersurfaces in toric variety",
math.DG/0007028.
%
\bibitem{SYZ}
Strominger, A.,Yau, S.-T. and Zaslow, E,
"Mirror Symmetry is T-duality",
{\it Nuclear Physics} B 479 (1996),243-259.
%
\bibitem{mirrorbook}
Yau, S.-T., (ed.), {\it Essays on Mirror Manifolds}, International Press, Hong
  Kong, 1992.
%
\bibitem{Z}
Zharkov, I.,
"Torus Fibrations of Calabi-Yau Hypersurfaces in Toric Varieties and Mirror Symmetry",
alg-geom 9806091
\end{thebibliography}
\end{document}